\documentclass[11pt,a4paper,reqno]{amsart}

\usepackage{lastpage}
\usepackage{amsmath}
\usepackage{enumerate}
\usepackage{euscript}
\usepackage{mathtools}
\usepackage{relsize}
\usepackage{latexsym}
\usepackage{amsfonts}
\usepackage{xcolor}
\usepackage{graphicx}
\usepackage{mathrsfs}
\usepackage{tikz}
\usepackage{accents}
\usepackage{pdflscape}
\usepackage{array}
\usepackage{booktabs}
\usepackage{multirow,rotating,caption}
\usepackage{appendix}
\usepackage{placeins}
\usepackage[numbers]{natbib}

\usepackage[letterpaper,top=4cm,bottom=4cm,left=3cm,right=3cm,marginparwidth=1.75cm]{geometry}
\usepackage[colorlinks=true, allcolors=blue]{hyperref}

\pgfdeclarelayer{edgelayer}
\pgfdeclarelayer{nodelayer}
\pgfsetlayers{edgelayer,nodelayer,main}

\tikzstyle{none}=[inner sep=0pt]

\newtheorem{thm}{Theorem}
\newtheorem{cor}[thm]{Corollary}
\newtheorem{lem}[thm]{Lemma}

\newtheorem{defn}[thm]{Definition}
\newtheorem{rem}[thm]{Remark}

\DeclareFontFamily{U}{mathx}{}
\DeclareFontShape{U}{mathx}{m}{n}{<-> mathx10}{}
\DeclareSymbolFont{mathx}{U}{mathx}{m}{n}
\DeclareMathAccent{\widehat}{0}{mathx}{"70}
\DeclareMathAccent{\widecheck}{0}{mathx}{"71}

\newlength{\widthofexpr}
\DeclarePairedDelimiterX{\setdef}[2]{\lbrace}{\rbrace}{#1\;\delimsize|\;\mathopen{}#2}

\renewcommand{\Re}{\operatorname{Re}}

\newcommand{\dom}{\operatorname{dom}}

\newcommand{\ran}{\operatorname{ran}}

\newcommand{\grad}{\nabla}
\newcommand{\divergence}{\operatorname{div}}
\newcommand{\curl}{\operatorname{curl}}

\newcommand{\R}{\ensuremath{\mathbb R}}
\newcommand{\C}{\ensuremath{\mathbb C}}

\newcommand{\N}{\ensuremath{\mathbb N}}

\newcommand{\calA}{\mathcal A}
\newcommand{\calB}{\mathcal B}

\newcommand{\calG}{\mathcal G}
\newcommand{\calH}{\mathcal H}

\newcommand{\calK}{\mathcal K}

\newcommand{\calR}{\mathcal R}

\newcommand{\calU}{\mathcal U}
\newcommand{\calV}{\mathcal V}

\newcommand{\calX}{\mathcal X}
\newcommand{\calY}{\mathcal Y}
\newcommand{\calZ}{\mathcal Z}

\newcommand{\spvek}[2]{\left(\begin{smallmatrix}#1\\#2\end{smallmatrix}\right)}
\newcommand{\sbmat}[4]{\left[\begin{smallmatrix}#1 & #2\\#3 & #4\end{smallmatrix}\right]}

\newcommand{\stack}[2]{\left[ \begin{smallmatrix} #1 \\ #2 \end{smallmatrix} \right]}

\newcommand{\roundstack}[2]{\left( \begin{smallmatrix} #1 \\ #2 \end{smallmatrix} \right)}

\newcommand{\nextto}[2]{\left[ \begin{smallmatrix} #1 \hspace{0.5mm} ,  \hspace{1mm} #2 \end{smallmatrix} \right]}

\newenvironment{smallbmatrix}{\left[\begin{smallmatrix}}{\end{smallmatrix}\right]}

\title{Abstract second-order boundary control systems}
\thanks{This is the accepted manuscript of an article published in Evolution Equations and Control Theory \textbf{27} (2027), 100835, \href{https://doi.org/10.3934/eect.2026.100835}{doi:10.3934/eect.2026.100835}.}

\author{Till Preuster$^1$}
\thanks{$^1$Junior Professorship Numerical Mathematics, Faculty of Mathematics, Chemnitz University of Technology, Germany, Mail: \textsc{\{till.preuster,manuel.schaller\}@math.tu-chemnitz.de}}
\author{Timo Reis$^2$}
\thanks{$^2$Fachgebiet Systemtheorie und partielle Differentialgleichungen, Institute of Mathematics, Technische Universität Ilmenau, Germany, Mail: \textsc{timo.reis@tu-ilmenau.de}}
\author{Manuel Schaller$^1$}

\thanks{This work was funded by the Deutsche Forschungsgemeinschaft (DFG, German Research Foundation) – Project-ID 531152215 – CRC 1701.}

\begin{document}

\maketitle

    \begin{abstract}
        Many physical systems, such as vibrating elastic bodies, acoustic media and electromagnetic fields, are governed by second-order evolution equations $\ddot x+D\dot x+Sx=0$ and interact with their environment through their boundary. The classical analysis via the operator matrix $\mathcal{A}=\left[\begin{smallmatrix}0&I\\-S&-D\end{smallmatrix}\right]$ requires $S$ to be self-adjoint, which typically fixes a boundary condition a priori and thus excludes boundary control. We therefore use $S$ only to define the energy space and govern the dynamics by a maximal operator $S_{\max}\supset S$ whose domain carries no boundary condition. Based on a factorization $S=A^*A$, we construct a boundary triplet for the associated operator matrix and define boundary inputs and outputs of Robin type, parametrized by an operator $P$ on the boundary space. We show that the resulting system is impedance passive and internally well-posed if $P$ is a contraction, and that for $D=0$ this covers all boundary conditions under which the free dynamics generate a contraction semigroup. Moreover, we show that the position-momentum formulation (e.g., $\left[\begin{smallmatrix}0&I\\\Delta&0\end{smallmatrix}\right]$ for the wave equation) and the strain-momentum formulation (e.g., $\left[\begin{smallmatrix}0&\nabla\\\operatorname{div}&0\end{smallmatrix}\right]$) describe the same open system, including its boundary inputs and outputs. We illustrate our results by an $n$-dimensional wave equation and Maxwell's equations.
    \end{abstract}

	\section{Introduction}
	Many physical systems obey the second-order in time evolution equation
	\begin{align}\label{eq:second_order_intro}
		\Ddot{x}(t)+ D \dot{x}(t) + Sx(t)=0, \qquad t \ge 0, \qquad x(0)=x_0, \quad \dot{x}(0)=x_1,
	\end{align}
	with initial position $x_0$ and initial velocity $x_1$, where the stiffness operator $S: \calX \supset \dom S \to \calX$ is self-adjoint and uniformly positive in a Hilbert space $\calX$ and $D\in L(\calX,\calX)$ is such that $-D$ is dissipative. Here, $S$ is a self-adjoint realization of a (formally symmetric) differential expression, that is, usually a boundary condition is encoded in $\dom S$. For instance, for the wave equation one may take the Dirichlet Laplacian $S=-\Delta_{\mathrm D}$ with $\dom S=\{x\in H_0^1(\Omega)\mid \Delta x\in L^2(\Omega)\}$, or for Maxwell's equations, $S$ is a self-adjoint realization of $c^2\curl\curl$ on divergence-free fields. The positivity of $S$ yields that $\calX_h = \dom S^{1/2}$ endowed with the energy norm $||S^{1/2} \cdot ||$ is a Hilbert space.
	Typically, the second-order evolution system \eqref{eq:second_order_intro} is analyzed via the associated first-order formulation
	\begin{align}\label{eq:calAdef}
		\calA = \begin{bmatrix}
			0 & I \\ -S & -D
		\end{bmatrix}, \qquad \dom \calA= \left\lbrace \roundstack{x}{y} \in \stack{\calX_h}{\calX_h} \, \middle| \, Sx+Dy \in \calX \right\rbrace
	\end{align}
	as an operator matrix in $\stack{\calX_h}{\calX}$. The operator $\calA$ is well-understood and research subject for more than 30 years, see \cite{banks1995well} and \cite{weiss2003get,jacob2008analyticity,jacob2017systems} concerning spectral properties of $\calA$. In particular, since $D$ is bounded, $\calA$ is a bounded perturbation of the skew-adjoint operator $\begin{smallbmatrix}
	    0 & I \\ -S & 0
	\end{smallbmatrix}$ and therefore even generates a $C_0$-group on $\stack{\calX_h}{\calX}$. Moreover, if additionally, $D$ is uniformly positive, then this group is exponentially stable. We emphasize that these facts relate to the closed system and are therefore not the subject of this article. However, the majority of works consider these systems closed, that is, without control input or observation and thus with a priori chosen boundary conditions.

	In this paper, we systematically associate boundary input and output maps to the underlying operator matrix $\calA$, thereby opening the previously closed system \eqref{eq:second_order_intro} towards interaction. Opening the system, e.g., to augment the system with a control, requires some care. Precisely, if $S$ is self-adjoint, a boundary condition is already contained in $\dom S$, and prescribing an additional boundary input could overdetermine the problem. Therefore, in the open system below, $S$ is replaced by a maximal operator $S_{\max}\supset S$ which acts by the same differential expression, but whose domain carries no boundary condition. The self-adjoint operator $S$ is only used to fix the energy space $\calX_h=\dom S^{1/2}$.
	Particularly, as a first contribution, we show well-posedness of
	\begin{align*}
		Gz(t)=u(t), \qquad \dot{z}(t)=\begin{bmatrix}
			0 & I \\ -S_{\max} & -D
		\end{bmatrix} z(t), \qquad K z(t)=y(t), \qquad z(0)=z_0
	\end{align*}
	for a predefined sufficiently regular control function $u$ with values in a boundary space and initial value $z_0$. The authors decided to work within the setting of boundary nodes provided in \cite{malinen2007impedance}. Therein, the operators $G$ and $K$ are for the moment not further specified input and output operators, typically boundary traces. To introduce them, we define a so-called boundary triplet $(\calG, \Gamma_0, \Gamma_1)$ for the operator matrix $\calA_\mathrm{max}=\sbmat{0}{I}{-S_\mathrm{max}}{0}$ where $\calG$ is the boundary space and $\Gamma_0, \Gamma_1$ map onto $\calG, \calG^*$, respectively. We then construct the boundary input and output mappings $G$ and $K$ based on linear combinations of $\Gamma_0$ and $\Gamma_1$. To illustrate this, consider the wave equation $\rho\ddot x+b\dot x-\divergence(T\grad x)+ax=0$ on a bounded Lipschitz domain $\Omega\subset\R^n$ with position $x$ and momentum $\rho\dot x$ as state variables. Here, the energy space is $H^1(\Omega)$, which carries no boundary condition, and the boundary maps are the boundary velocity $\Gamma_0 z=\dot x|_{\partial\Omega}$ and the Neumann trace $\Gamma_1 z=\nu\cdot T\grad x|_{\partial\Omega}$.
    The inputs considered in this article are of the form $u=\tfrac12(1-P)\calR\,\dot x|_{\partial\Omega}+\tfrac12(1+P)\,\nu\cdot T\grad x|_{\partial\Omega}$, where $P$ is a contraction on the boundary space and $\calR$ is the Riesz isomorphism on $H^{1/2}(\partial\Omega)$. They comprise Neumann boundary control for $P=1$ and control of the boundary velocity for $P=-1$. For a suitable choice of $P$, they also comprise the classical dissipative boundary control $\nu\cdot T\grad x+k\dot x=\tilde u$ (on parts of the boundary) with $k\ge0$, which is used for boundary stabilization of the wave equation \cite{lagnese1983decay,komornik1990direct,lasiecka2000control}. Having defined such a boundary triplet for $\calA_\mathrm{max}$, all possible configurations of boundary nodes leading to well-posed impedance passive boundary control systems are derived.
	We stress that there is a rich theory concerning control and observation of (conservative) physical systems on infinite-dimensional state spaces. In particular, the interplay between boundary triplets and boundary control systems has been thoroughly investigated; see \cite{arov2012boundary} and for the connection between state/signal theory and boundary relations we recommend \cite{arov2012passive}.
	In the present contribution, we adopt and extend the techniques developed in \cite{malinen2007impedance,malinen2006conservative} to the setting of second-order systems. In particular, we introduce a new construction showing that boundary triplets for the governing operators of abstract second order systems arise as transformations of boundary triples for so-called dual pairs, a relation that has not been developed in \cite{malinen2007impedance,wegner2017boundary}. Moreover, we note that the framework of system nodes, see \cite{Staffans2005,PhilReis23}, allows for a unified treatment of combined distributed and boundary control.
        Boundary control of second-order hyperbolic equations has been developed to great depth by PDE methods, including sharp regularity for Dirichlet and Neumann boundary data \cite{lasiecka1986non,lasiecka1991regularity,tataru1998regularity}, exact controllability \cite{lions1988exact}, uniform boundary stabilization \cite{lagnese1983decay,komornik1990direct}, and optimal control and Riccati theory, see \cite{lasiecka2000control} and the references therein. These results rely on specifically PDE-based techniques and are, in general, not accessible by abstract semigroup theory alone - the main subject of this work.

	As a second contribution, we provide equivalence transforms of the above system class which in particular appear in port-Hamiltonian formulations.
	Port-Hamiltonian systems provide an modular energy-based modeling approach of dissipative systems that are subject to impedance passive control and observation. For works on linear infinite-dimensional formulations, we refer to \cite{JacoZwar12,villegas2007port} for systems on one-dimensional spatial domains, and to~\cite{skrepek2021linear} and \cite{PhilReis23} for systems on multi-dimensional domains. In port-Hamiltonian formulations, instead of using the standard first-order reformulation $\calA$ of \eqref{eq:second_order_intro}, a reformulation building upon an factorization $S=A^*A$ is often used.
	Introducing the state transformation $z= \stack{Ax}{\dot{x}}$ one obtains, in contrast to $\calA$ as defined in \eqref{eq:calAdef}, the operator matrix
	\begin{align*}
		\calB = \begin{bmatrix}
			0 & A \\ -A^* & -D
		\end{bmatrix}.
	\end{align*}
	To illustrate this, in case of the wave equation one has $S=-\Delta=-\divergence ( \grad \cdot )$. To show well-posedness and to derive an associated boundary control system, the authors in \cite{kurula2015linear} consider the generator $\calB=\sbmat{0}{\grad}{\operatorname{div}}{0}$. In this regard we also mention the recent contribution \cite{aigner2025well} achieving a similar goal choosing the generator $\calA$ from \eqref{eq:calAdef} instead. We prove that the two formulations $\calA$ and $\calB$ are related by an energy-induced similarity in general. On a formal level, this type of system transformation has been studied in \cite{preuster2024jet} for port-Hamiltonian partial differential equations on one-dimensional spatial domains. Furthermore, at the level of boundary triplets, there is a one-to-one correspondence between the maximal operators associated with $\calA$ and $\calB$ which paves the way for a boundary control system associated with $\calB$. Thusly, the similarity of the two operator matrices bridges the gap between the two formulations and enables the transfer of properties such as well-posedness, passivity, and energy balance from one framework to the other.

    \medskip
	This paper is organized as follows. In Section \ref{sec:prel} we recall the necessary knowledge from the mathematical domains used in this paper. In particular, we introduce basic notions from unbounded operators and forms. Thereafter, existing theory of boundary nodes and boundary triplets is recalled. We present the main contributions in Section \ref{sec:main}. In Subsection~\ref{subsec:abpot}, we derive the abstract setting of the underlying state spaces used in the formulation of second-order systems. Then, in Subsection~\ref{subsec:bdd_triplet}, we explicitly construct a boundary triplet for the governing operator $\calA_\mathrm{max}$. A canonical equivalence transform in the sense of the factorization of $S$ is given in Subsection~\ref{subsec:Helmholtz}. We then provide the definition of the boundary node associated with the second-order system \eqref{eq:second_order_intro} in Subsection~\ref{subsec:bdd_node} and we allow for the presence of nontrivial mass-density and dissipation operators in Subsection~\ref{subsec:Mda}. The application of the theory is provided by means of a wave and  Maxwell equation in Section \ref{sec:applications}. We conclude our work in Section~\ref{sec:conclusion}. Well-known trace theorems for standard differential  operators in Sobolev spaces are provided in Appendix~\ref{sec:app}. To make the novelty transparent, we note that Section~\ref{sec:prel} only collects known material, whereas Section~\ref{sec:main} and the applications in Section~\ref{sec:applications} are new.

	\section{Preliminaries}\label{sec:prel}

	In this section, we recall the fundamental concepts that will be used throughout the paper. All results of this section are known and are included to keep the article self-contained. Specifically, Subsection~\ref{subsec:operators_form_rel} reviews essential notions from functional analysis, Subsection~\ref{subsec:sysnodes} introduces the framework of boundary nodes,
	and Subsection~\ref{subsec:prel_bdd_triplets} summarizes key ideas from the theory of boundary triplets.
	\subsection{Operators, forms and relations in Hilbert spaces}\label{subsec:operators_form_rel}
	All Hilbert spaces are complex throughout this article. For a~Hilbert space $\calX$, the {\em anti-dual}, i.e., the set of bounded and conjugate linear functionals on $\calX$, is denoted by $\calX^*$. If not stated otherwise, the anti-dual is identified with $\calX$, which is possible by the Riesz representation theorem
	\cite[Thm.~6.1]{alt2016linear}. The inner product and norm in $\calX$ are denoted by $\langle\cdot,\cdot\rangle_\calX$ and $\|\cdot\|_\calX$, respectively, whereas
	$\langle\cdot,\cdot\rangle_{\calX^*,\calX}$ is the duality pairing between $\calX$ and $\calX^*$.
	The (linear) Riesz isomorphism from $\calX$ to $\calX^*$
	is denoted by $\calR_{\calX}$. Throughout this article, the Cartesian product of two Hilbert spaces  $\calX$ and $\calY$ is denoted by $\stack{\calX}{\calY}$ and is endowed with the usual product Hilbert space norm.

	An alternative consideration of anti-duals of Hilbert spaces arises
	naturally in the framework of \emph{Gelfand triples}.
	A Gelfand triple consists of a chain of continuous embeddings
	\[
	\calV \hookrightarrow \calX \hookrightarrow \calV^*,
	\]
	where $\calV$ is another Hilbert space densely and continuously embedded in  $\calX$.
	Here, $\calX$ is called the \emph{pivot space}, since it is identified with its own anti-dual.
	This allows us to uniquely extend the restriction of $\langle \cdot,\cdot\rangle_{\calX}$ to $\stack{\calV}{\calX}$ into a duality pairing between $\calV$ and $\calV^*$.
	For further details on Gelfand triples, we refer to the original paper \cite{gelfand_kostyucenko_1955}.

	For a~further Hilbert space $\calY$, the space of bounded linear operators from $\calX$ to $\calY$ is denoted by $L(\calX,\calY)$, and we abbreviate $L(\calX)\coloneqq  L(\calX,\calX)$. Moreover, $I_\calX$ stands for the identity on $\calX$.

	A (possibly unbounded) operator $A: \calX \supset \dom A \to \calY$ is {\em closed}, if the {\em domain} $\dom A$ is complete with respect to the {\em graph norm}
	\[\|x\|_{\dom A}\coloneq \big(\|x\|_\calX^2+\|Ax\|_\calY^2\big)^{1/2}.\]
	For operators $A: \calX \supset \dom A \to \calY$, $B: \calX \supset \dom B \to \calY$, we write $A\subset B$, if $\dom A\subset\dom B$, and if $A$ is the restriction of $B$ to $\dom A$.\\
	The \textit{adjoint} $A^*:\calY^* \supset\dom A^* \to\calX^*$ of a~densely defined linear operator $A:\calX\supset\dom A\to\calY$ has the domain
	\[
	\dom A^*=\setdef[\big]{y^*\in \calY^*}{\exists\, z^*\in\calX^* \text{ s.t.\ }\forall\,x\in\dom A:\;\langle y^*,Ax\rangle_{\calY^*,\calY}=\langle z^*,x\rangle_{\calX^*,\calX}}.
	\]
	The density of $\dom A$ in $\calX$ implies that $z$ in the above set is uniquely determined by $x$, which ensures that $A^*y = z$ is well-defined.
	An operator $A \colon  \calX\supset\dom A \to \calX$ is called \emph{symmetric} if $A \subset A^*$ and \emph{self-adjoint} if $A = A^*$.
	Analogously, $A$ is called \emph{skew-symmetric} if $A \subseteq -A^*$ and \emph{skew-adjoint} if $A = -A^*$.

	A mapping $h \colon \dom h \times \dom h \to \C$, where $\dom h$ is a subspace of $\calX$, is called a \emph{sesquilinear form} if it is linear in its first argument and conjugate-linear in its second.
	The form $h$ is called \emph{symmetric} if $h(x,y) = \overline{h(y,x)}$ for all $x,y \in \dom h$.
	A symmetric form is called \emph{positive} if $h(x,x) > 0$ for all $x \in \dom h \setminus \{0\}$, and \emph{coercive} if there exists a constant $c>0$ such that $h(x,x)>c\,\|x\|^2_\calX$ for all $x\in\dom h$. We shortly write $h(x)\coloneqq  h(x,x)$ for $x \in \dom h$.
	For a~positive form $h$, the mapping
	\[x\mapsto \|x\|_{\dom h} \coloneq \big(\|x\|^2_\calX+h(x)\big)^{1/2}\]
	is a~norm on $\dom h$ (which obviously originates from an inner product on $\dom h$).
	If $\dom h$ is complete with respect to this norm, then we call $h$ {\em closed}.

	For a self-adjoint operator $P \colon \calX \supset \dom P \to \calX$, the mapping $x \mapsto \langle x, Px \rangle_\calX$ defines a symmetric form.
	We say that $P$ is \emph{positive} (respectively, \emph{uniformly positive}) if the associated form is positive (respectively, coercive). Any positive $P \colon \calX \supset \dom P \to \calX$ possesses a~uniquely determined {\em operator square root} $P^{1/2} \colon \calX \supset \dom P^{1/2} \to \calX$, that is, $P^{1/2}$ is positive with $(P^{1/2})^2=P$.  Thus, the form
	$x \mapsto \langle x, Px \rangle_\calX$ extends to a~closed form $h$ with $\dom h=\dom P^{1/2}$ and $h(x)=\|P^{1/2}x\|_\calX^{2}$. For further details on forms, we refer to \cite[Chap.~6]{kato2013perturbation}.

	A \emph{linear relation} $\calK$ on $\calX$ is a subspace of $\stack{\calX}{\calX}$.
	The relation $\calK$ is called \emph{dissipative} if
	\[
	\Re \langle f, g \rangle_{\calX} \leq 0
	\quad \text{for all } \spvek{f}{g} \in \calK.
	\]
	If $\calK$ admits no proper dissipative extension, it is called \emph{maximal dissipative}.

	An operator $A \colon \dom A \subset \calX \to \calX$ is called (maximal) dissipative
	if its \emph{graph}
	$\setdef{({x},{Ax})}{ x \in \dom A }$
	is (maximal) dissipative.

	Throughout this article we follow the conventions of \cite{adams2003sobolev} concerning the notation of Lebesgue and Sobolev spaces.
	When dealing with function spaces whose values lie in a Hilbert space $\calX$, we indicate the target space by appending “;$\calX$” after the domain.
	For instance, the space of $\calX$-valued functions on $\Omega \subset \R^d$ that are $p$-integrable is written as $L^{p}(\Omega;\calX)$.
	In this article, all integrals of $\calX$-valued functions are to be understood in the sense of Bochner (see \cite{Dies77}).

	\subsection{Boundary nodes}\label{subsec:sysnodes}
	To introduce an abstract framework for boundary control systems, we adopt the notion of \textit{boundary nodes} as developed by {\sc Malinen and Staffans} in \cite{malinen2007impedance}.
	The definitions and results presented in this subsection are taken from that reference and are included here for completeness. For Hilbert spaces $\calX$, $\calU$, and $\calY$, we consider systems of the form
	\begin{equation}\label{eq:ODEnode}
		u(t)=Gz(t), \qquad \dot{z}(t) =Lz(t), \qquad y(t)=Kz(t), \qquad z(0)=z_0
	\end{equation}
	where $L$ is typically some differential operator and $G$ and $K$ are boundary traces. Systems of this type naturally arise in the control theory of distributed-parameter and physical systems, see \cite{reis_infinite-dimensional_dissipation_2025,PhilReis23}. The concept of boundary nodes provides an operator-theoretic foundation for such systems and facilitates their treatment in real-world applications, cf. \cite{ReisSchaller2024Oseen}.
	\begin{defn}
		A \textit{colligation} $\Theta = \left(\begin{bsmallmatrix}
			G \\ L \\ K
		\end{bsmallmatrix}, \begin{bsmallmatrix}
			\calU \\  \calX  \\ \calY
		\end{bsmallmatrix}\right)$ consists of three Hilbert spaces $\calU, \calX$, and $\calY$, and three linear maps $G,L$, and $K$, with the same domain $\calZ \subset \calX$ and with values in  $\calU, \calX$, and $\calY$, respectively. The domain $\dom \Theta$ of $\Theta$ denotes the common domain $\calZ$ of $G,L$, and $K$. The colligation $\Theta$ is \textit{closed} if the combined operator
		\begin{align*}
			\begin{bsmallmatrix}
				G \\ L \\ K
			\end{bsmallmatrix} : \calX \supset \calZ \to \begin{bsmallmatrix}
				\calU \\  \calX  \\ \calY
			\end{bsmallmatrix}
		\end{align*}
		is closed.
	\end{defn}
	\begin{defn}[Boundary node]\label{def:boundarynode}
		A colligation $\Theta = \left(\begin{bsmallmatrix}
			G \\ L \\ K
		\end{bsmallmatrix}, \begin{bsmallmatrix}
			\calU \\  \calX  \\ \calY
		\end{bsmallmatrix}\right)$  is a {\em boundary node} if the following conditions are satisfied:
		\begin{enumerate}[(a)]
			\item $\Theta$ is closed;
			\item\label{def:boundarynode2} the operator $G$ is surjective and $\ker G$ is dense in $\calX$;
			\item the restriction $A \coloneqq  \left. L \right|_{\ker G}$ has non-empty resolvent set.
		\end{enumerate}
		If additionally,
		\begin{enumerate}
			\item[(d)]  the operator $A$ is the generator of a~strongly continuous semigroup $\mathfrak{A}(\cdot)\colon
			\R_{\ge 0}\to L(\calX)$ on $\calX$,
		\end{enumerate}
		this boundary node is \textit{internally well-posed} (in the forward time direction). We call $\calU$ the \textit{input space}, $\calX$ the \textit{state space}, $\calY$ the \textit{output space}, $\calZ$ the \textit{solution space}, $G$ the \textit{input boundary operator},  $L$ the \textit{interior operator}, $K$ the \textit{output boundary operator}, and $A$ the \textit{main operator} (or the \textit{semigroup generator} in the internally well-posed case).
	\end{defn}
	We recall the concept of classical solutions from \cite[Lem. 2.6]{malinen2006conservative}.
	\begin{lem}\label{lem:solex}
		Let $\Theta  = \left(\begin{bsmallmatrix}
			G \\ L \\ K
		\end{bsmallmatrix}, \begin{bsmallmatrix}
			\calU \\  \calX  \\ \calY
		\end{bsmallmatrix}\right)$ be an internally well-posed boundary node. Then, for all $z_0 \in \calX$ and $u \in C^2(\R_{\geq 0}; \calU)$ with $Gz_0=u(0)$ the first, second and fourth equation in \eqref{eq:ODEnode} have a unique solution $z \in C^1(\R_{\geq 0}; \calX) \cap C(\R_{\geq 0}; \calZ)$. Hence, we can define $y \in C(\R_{\geq 0}; \calY)$ by the third equation in \eqref{eq:ODEnode}. In this work, we shortly say “a smooth solution of \eqref{eq:ODEnode} on $\R_+$” when we mean a solution with the above properties.
	\end{lem}
	\begin{defn}[Scattering passive boundary node]
		A boundary node $\Theta$ on $(\calU,\calX, \calY)$ is \textit{scattering passive} if it is internally well-posed and all smooth solutions of \eqref{eq:ODEnode} on $\R_+$ satisfy
		\begin{align}\label{eq:passivity_ineq}
			\frac{\mathrm{d}}{\mathrm{d}t} ||z(t)||_\calX^2 +||y(t)||_\calY^2 \leq ||u(t)||_\calU^2, \qquad t \in \R_+.
		\end{align}
		It is \textit{scattering energy-preserving} if the above inequality holds in the form of an equality.
	\end{defn}
	\begin{defn}[External Cayley transform and impedance passive boundary node]\label{def:imped_bn}
		By the \textit{external Cayley transform the colligation $\Theta  = \left(\begin{bsmallmatrix}
				G \\ L \\ K
			\end{bsmallmatrix}, \begin{bsmallmatrix}
				\calU \\  \calX  \\ \calU
			\end{bsmallmatrix}\right)$ with parameter $\beta \in \C_+$}  we mean the colligation $\Theta^{(\beta)}  = \left(\begin{bsmallmatrix}
			G^{(\beta)} \\ L \\ K^{(\beta)}
		\end{bsmallmatrix}, \begin{bsmallmatrix}
			\calU  \vphantom{ G^{(\beta)}}\\  \calX \vphantom{L}  \\ \calU \vphantom{K^{(\beta)}}
		\end{bsmallmatrix}\right)$ with $\dom \Theta^{(\beta)} = \dom \Theta$ and
		\begin{align}
			G^{(\beta)}=\frac{1}{\sqrt{2\Re \beta}}(\beta G+K), \qquad \text{and} \qquad K^{(\beta)}=\frac{1}{\sqrt{2\Re \beta}}(\overline{\beta} G-K).
		\end{align}
		Let $\Theta  = \left(\begin{bsmallmatrix}
			G \\ L \\ K
		\end{bsmallmatrix}, \begin{bsmallmatrix}
			\calU \\  \calX  \\ \calU
		\end{bsmallmatrix}\right)$ be a colligation and $\Theta^{(\beta)}  = \left(\begin{bsmallmatrix}
			G^{(\beta)} \\ L \\ K^{(\beta)}
		\end{bsmallmatrix}, \begin{bsmallmatrix}
			\calU  \vphantom{ G^{(\beta)}}\\  \calX \vphantom{L}  \\ \calU \vphantom{K^{(\beta)}}
		\end{bsmallmatrix}\right)$ be the external Cayley transform of $\Theta$ with parameter $\beta \in \C_+$.
		\begin{enumerate}
			\item[(i)] $\Theta$ is \textit{impedance passive} if $\Theta^{(\beta)}$ is a scattering passive boundary node for some $\beta \in \C_+$.
			\item[(ii)] $\Theta$ is \textit{impedance energy-preserving} if $\Theta^{(\beta)}$ is a scattering energy-preserving boundary node for some $\beta \in \C_+$.
		\end{enumerate}
	\end{defn}
	We note that \cite[Proposition 4.1]{malinen2007impedance} provides a characterization of internal well-posedness of impedance passive colligations.

	\subsection{Boundary triplets and dual pairs}\label{subsec:prel_bdd_triplets}
	A suitable framework for abstractly formulating boundary conditions of partial differential equations is provided by the \emph{boundary triplets}, rooted in abstract extension theory of symmetric operators \cite{gorbachuk1989extension,derkach1991generalized}.
	We now introduce this approach and recall its most important consequences.\begin{defn}\label{def:operator_bdd_triplet}
		Let $\calX$ be a Hilbert space and $A_0 \colon \calX \supset \dom A \to \calX$ be a closed, densely defined and skew-symmetric operator. A \emph{boundary triplet} for $A_0^*$ is a triplet $\left( \calG, \Gamma_0, \Gamma_1\right)$ consisting of a Hilbert space $\calG$ and a linear surjective map $\Gamma =  \left[ \begin{smallmatrix}
			\Gamma_0 \\ \Gamma_1
		\end{smallmatrix} \right]: \dom A_0^* \to \stack{\calG}{\calG^*}$ such that the \textit{abstract Green identity}
		\begin{align}\label{eq:operator_Green}
			\left\langle A_0^*f , g \right\rangle_\calX + \left\langle f , A_0^*g \right\rangle_\calX = \left\langle \Gamma_1 f, \Gamma_0 g\right\rangle_{\calG^*,\calG} + \left\langle \Gamma_0f, \Gamma_1 g\right\rangle_{\calG,\calG^*}
		\end{align}
		holds for all $f,g \in \dom A_0^*$.
	\end{defn}
	Boundary triplets are typically formulated in the setting where $\calG$ is identified with its anti-dual.
	However, by the applicability of the Riesz isomorphism, the theory and results remain exactly the same in the case where $\calG^*$ is not identified with $\calG$.

	Let $(\calG, \Gamma_0, \Gamma_1)$ be a boundary triplet for an operator $A_0^*$.
	It can be verified that $\Gamma \colon \dom A_0^* \to \stack{\calG}{\calG^*}$ is bounded,
	where $\dom A_0^*$ is equipped with the graph norm of $A_0^*$, \cite[Prop. 2.1.2]{BehrHdS20}.

	Also the domain of the \textit{minimal operator} $A_0$ can be recovered by restricting the domain of the \textit{maximal operator} $A_0^*$ to the kernel of the boundary map $\Gamma$.
	\begin{lem}[{\cite[Lem.\ 4.3]{wegner2017boundary}}]\label{lem:cont_bdd_map}
		Let $A_0 \colon \calX \supset \dom A \to \calX$ be a closed, densely defined and skew-symmetric operator. Moreover, let $\left( \calG, \Gamma_0, \Gamma_1\right)$ be a boundary triplet for $A_0^*$. Then
		\[
		\dom A_0=\ker \Gamma = \ker \Gamma_0 \cap \ker \Gamma_1.
		\]
	\end{lem}
	It turns out that there is a one-to-one correspondence between the maximal dissipative extensions of $A_0$ and contractions $P$ on the boundary space $\calG^*$.
	\begin{thm}[{\cite[Thm.~4.2]{wegner2017boundary}}]\label{thm:former_Kurula_Zwart}
		Let $A_0 \colon \calX \supset \dom A_0 \to \calX$ be a closed, densely defined and skew-symmetric operator. Let $\left( \calG, \Gamma_0, \Gamma_1\right)$ be a boundary triplet for the operator $A_0^*$. Then, an extension $B: \calX \supset \dom B \to \calX$ of $A_0$ generates a contractive semigroup on $\calX$ if and only if there exists a uniquely determined $P \in L(\calG^*)$ with $\|P\|_{L(\calG^*)} \leq 1$ such that
		\begin{align*}
			\dom B =  \left\lbrace x \in\dom A_0^* \, \middle| \,(P-1)\calR_{\calG} \Gamma_0 x - (P+1) \Gamma_1  x = 0 \right\rbrace
		\end{align*}
		and $B=\left. A_0^* \right|_{\dom B}$.
	\end{thm}
	The following definition is slightly adapted from \cite{malamud2002krein}, see also \cite{derkach2023pt}.
	\begin{defn}[Dual pairs]\label{def:dp}
		Let $\calX, \calY$ be Hilbert spaces and let $M: \calY \supset \dom M \to \calX$ and $N: \calX \supset \dom N \to \calY$ be closed and densely defined linear operators. Then, $(M,N)$ is a \textit{dual pair} if
		\begin{align*}
			\forall\,y \in \dom M,\, x \in \dom N:\quad    \langle My,x \rangle_\calX - \langle y, Nx \rangle_\calY =0
		\end{align*}
		Let $\calG_1, \calG_2$ be Hilbert spaces. Let
		\begin{align}
			\Lambda = \stack{\Lambda_1}{\Lambda_2}: \dom M^* \to \stack{\calG_1}{\calG_2^*}, \quad \Pi = \stack{\Pi_1}{\Pi_2}: \dom N^* \to \stack{\calG_1^*}{\calG_2}
		\end{align}
		be linear operators. The triple $(\stack{\calG_1}{\calG_2}, \Lambda, \Pi)$ is called a \textit{boundary triplet for the dual pair $(M,N)$} if
		\begin{enumerate}
			\item[i)] $\ran \Lambda =  \stack{\calG_1}{\calG_2^*}$, $\ran \Pi = \stack{\calG_1^*}{\calG_2}$;
			\item[ii)] the abstract Green identity \begin{align}
				\langle N^*y,x \rangle_\calX - \langle y, M^*x \rangle_\calY = \langle \Pi_1 y, \Lambda_1 x \rangle_{\calG_1^*, \calG_1} - \langle \Pi_2 y, \Lambda_2 x \rangle_{\calG_2, \calG_2^*}
			\end{align}
			holds for all $y \in \dom N^*, x \in \dom M^*$.
		\end{enumerate}
	\end{defn}
	Note that $(M,N)$ being a~dual pair is equivalent to $M \subset N^*$, which is in turn equivalent to $N \subset M^*$.

	Moreover, as we will see later in Theorem~\ref{thm:bdd_triplet_for_sec_order}, there is a closed connection between boundary triplets for dual pairs (as introduced in above Definition~\ref{def:dp}) and boundary triplets for a suitably chosen block operator $\calA_\mathrm{max}$ (as defined in Definition~\ref{def:operator_bdd_triplet}).

	\section{Abstract second order systems}\label{sec:main}
	Let $\calX$ be a Hilbert space. For motivational purposes, we start by considering the free dynamics (i.e., without control input) of second-order systems, namely
	\begin{align}\label{eq:second_order}
		M\Ddot{x}(t)+ D \dot{x}(t) + Sx(t)=0
	\end{align}
	where $S: \calX \supset \dom S \to \calX$ and $M\in L(\calX)$ are self-adjoint and uniformly positive, whereas $D\in L(\calX)$ with $-D$ dissipative. As explained in the introduction, $S$ defines a fixed self-adjoint realization and only serves to define the underlying energy space. Define the Hilbert space
	\begin{align}\label{eq:defX_h}
		\calX_h=\dom S^{1/2}\  \text{endowed with the norm}\ \|S^{1/2} \cdot \|_\calX
	\end{align}
	and let $\calX_h^*$ be the dual space of $\calX_h$ with respect to the pivot space $\calX$. Positivity of $S$ yields that the norm on $\calX_h$ is stronger than the norm in $\calX$.
	Further, $S$ restricts to a bounded linear operator $S \colon \calX_h \to \calX_h^*$, which coincides with the Riesz isomorphism on $\calX_h$. By boundedness and positivity, the inner product $\langle M^{-1} \cdot ,\cdot\rangle_\calX$ is equivalent to the standard inner product on $\calX$. When writing $\calX_M$ we mean the space $\calX$ equipped the above inner product $\langle M^{-1} \cdot ,\cdot\rangle_\calX$.

	We introduce the state variables
    \begin{align}\label{eq:Z_h}
        z(t)=\spvek{x(t)}{M\dot{x}(t)}\in  \stack{\calX_h }{ \calX_M} \eqqcolon \calZ_h
    \end{align}
	which may be viewed as abstract position and momentum.
	Then, we rewrite \eqref{eq:second_order} as an abstract Cauchy problem, \begin{align}\label{eq:Cauchy_second_order}
		\tfrac{\mathrm{d}}{\mathrm{d}t} z(t) = \calA z(t), \qquad z(0)=z_0
	\end{align}
	where $z_0$ contains the initial values for $x$ and $M\dot{x}$, and $\calA \colon \calZ_h \to \calZ_h$ with
	\begin{align}\label{eq:As}
		\dom \calA = \setdef{ \spvek{z_1}{z_2} \in   \stack{\calX_h }{ \calX_h}}{ S z_1 \in \calX_M}, \qquad \calA z= \sbmat{
			0}{M^{-1}}{-S}{-DM^{-1}}
		z,
	\end{align}
	where we note that in the definition of $\dom \calA$, $S$ has to be understood as an operator from $\calX_h$ to $\calX_h^*$ in the sense of the above discussion after \eqref{eq:second_order}.
	Operators of the form $\calA$ are well-understood, see \cite{weiss2003get}, and related analyses also exist for the case of unbounded damping operators, see, for instance, \cite{jacob2008analyticity,jacob2017systems}. In particular,  such operators are maximal dissipative and they generate a contraction semigroup on $\calZ_h$. If $D=0$, they are skew-adjoint and they generate a~unitary semigroup on $\calZ_h$.

	Motivated by second-order equations from mechanics, we introduce the {\em abstract potential energy} $H_p:\calX_h\to\R$, the {\em abstract kinetic energy} $H_k:\calX_M\to\R$ and {\em abstract total energy} $H:\stack{\calX_h }{ \calX_M}\to\R$
	of \eqref{eq:Cauchy_second_order} by
	\begin{align}\label{eq:energies}
		H_p(z_1)=\tfrac{1}{2} \|z_1\|_{\calX_h}^2,\quad  H_k(z_2)= \tfrac{1}{2} \|z_2\|_{\calX_M}^2 = \tfrac{1}{2} \langle M^{-1} z_2,z_2 \rangle_\calX,\quad
		H(z)= H_p(z_1) + H_k(z_2).
	\end{align}
	Note that the operator
	\[
	\calH \colon  \calZ_h \to \calZ_h^*=\stack{\calX_h^*}{ \calX_M},
	\qquad
	\calH = \begin{bmatrix} S & 0 \\[1mm] 0 & M^{-1} \end{bmatrix},
	\]
	is the Gâteaux derivative of the energy functional $H$ with respect to $z$.
	In this sense, the operator
	\[
	\calA = \left(
	\begin{bmatrix} 0 & I \\[1mm] -I & 0 \end{bmatrix} -
	\begin{bmatrix} 0 & 0 \\[1mm] 0 & D \end{bmatrix} \right)
	\begin{bmatrix} S & 0 \\[1mm] 0 & M^{-1} \end{bmatrix}
	\]
	represents an abstract (dissipative) Hamiltonian law, i.e., the derivative of the energy is left-multiplied with the canonical skew-symmetric operator.
	Interpreting the system~\eqref{eq:second_order} in terms of abstract energies
	is fundamental for defining the underlying state space of~\eqref{eq:Cauchy_second_order}.
	\par
	Next, we formulate abstract and suitable boundary conditions. Hereby, we will drop the assumption of $S$ being self-adjoint to allow the maximally possible freedom in terms of boundary conditions. We will merely assume that $S$ possesses a~self-adjoint and coercive restriction. As we will see, this allows us to parametrize all suitable boundary input-output configurations
	that yield an impedance passive boundary control system (i.e., a port-Hamiltonian system).

	To this end, we first investigate the parametrization of the abstract potential energy in Subsection~\ref{subsec:abpot}. While assuming that $M=I$ and $D=0$ we derive a boundary triplet for the maximal operator associated to \eqref{eq:second_order} in Subsection~\ref{subsec:bdd_triplet}. Afterwards, in Subsection~\ref{subsec:Helmholtz}, we apply a canonical system transformation, based on a factorization of $S$, while preserving its physical interpretation. Then, we derive the associated boundary node in Subsection~\ref{subsec:bdd_node}.
	Finally, in Subsection~\ref{subsec:Mda}, we will allow for the presence of nonzero damping and more general (abstract) mass density operators $M$.

	\subsection{Abstract potential energy and associated Gelfand triple}\label{subsec:abpot}
	Our main motivation comes from concrete examples in partial differential equations. In many such cases, after choosing suitable boundary conditions, $S$ is a self-adjoint uniformly positive operator $S$ and admits a factorization $S=A^*A$, where $A$ is a differential operator that maps into a second Hilbert space $\calY$. Typical examples include $A=\grad$ with $\dom A=H^1(\Omega)$, leading to the Neumann Laplacian (augmented by a zero-order term to ensure coercivity, see Section~\ref{subsec:rod}), $A=\grad$ with $\dom A=H_0^1(\Omega)$, leading to the Dirichlet Laplacian, and $A=c\curl$ on suitable spaces of divergence-free fields, leading to a self-adjoint realization of $c^2\curl\curl$ (Maxwell's equations, see Section~\ref{subsec:maxwell}). Note that the boundary condition contained in $\dom S$ for $S=A^*A$ is determined by the choice of $\dom A$. As observed at the beginning of this section, the form $h$ induces a state space for the position variable $x$. We shall subsequently study a system transformation whose underlying space is $\stack{\calY}{\calX}$, with transformation map defined by the operator $\sbmat{A}{}{}{1}$. For example, the similarity between the generator of the wave equation in \emph{first-order form} $\sbmat{0}{\grad}{\operatorname{div}}{0} $ and the \textit{second-order} generator $\sbmat{0}{I}{\Delta}{0}$ is formally given by $\sbmat{A}{}{}{1}=\sbmat{\grad}{}{}{1}$. Precisely, the state transformation $y=\sbmat{\grad}{}{}{1}z$ yields the governing operator
	\begin{align}\label{eq:wavebeispiel}
		\sbmat{\grad}{}{}{1} \sbmat{0}{I}{\Delta}{0}\sbmat{\grad}{}{}{1}^{-1} =\sbmat{0}{\grad}{\operatorname{div}}{0}
	\end{align}
	which corresponds to the transition from position-momentum to strain-momentum formulation.
	In what follows, we assume that there is a closed and densely defined linear operator $A: \calX \supset \dom A \to \calY$ such that the operator $S$ as in \eqref{eq:Cauchy_second_order} is factorized as $S=A^*A$. Note that such a factorization always exists in view of the operator square root, but may be chosen also differently. To enable more freedom in terms of boundary conditions and in particular for boundary control access, our first main objective is to give an extension $-B \supset A^*$ in order to replace $S$ by $-B A$. In particular, $(A^*,-B^*)$ will give rise to a dual pair from which we may derive a boundary triplet and an associated boundary control system.

	First, positivity of $S$ implies that $A$ is bounded from below by a strictly positive constant, i.e. there exists a constant $c>0$ such that $\|Ax\|_\calY \geq c \|x\|_\calX$ for all $x \in \dom A$. This is due to the fact that the operator square root $|A| \coloneqq  (A^*A)^{1/2}=S^{1/2}$ satisfies $\dom |A| = \dom A$ and $||\, |A|x||_\calX=||Ax||_{\calY}$ for all $x \in \dom A$, i.e. $|A|$ and $A$ are \textit{metrically equivalent}.
    Consequently, the Hilbert space $\calX_h$ coincides with $\dom A$ and is endowed with the inner product $h: \dom h \times \dom h \to \C$ defined by
	\begin{align}\label{eq:h}
		h(x,y) = \langle Ax,Ay \rangle_\calY , \qquad \dom h = \dom A.
	\end{align}
	Furthermore, the coercivity of the form $h$ directly yields that $A$ has closed range, see \cite[Thm.~IV.5.2]{kato2013perturbation}. This gives rise to the orthogonal decomposition
	\begin{align}\label{eq:Y_1}
		\calY = \calY_0 \oplus \calY_1 \eqqcolon (\ran A )^\perp \oplus \overline{\ran A} = \ker A^* \oplus \ran A.
	\end{align}
	Observe that $(\calY_1, \langle \cdot, \cdot \rangle_{\calY_1}) \coloneqq  (\ran A, \langle \cdot, \cdot \rangle_\calY)$ is a Hilbert space by closedness of $\ran A$ in $\calY$. Hence, $A$ restricts to a bounded linear operator, which we denote again by
	\begin{align}\label{eq:bounded_restriction}
		A: \calX_h \to \calY_1.
	\end{align}
	Moreover, $A$ as in \eqref{eq:bounded_restriction} is an isometric isomorphism from $\calX_h$ to $\calY_1$ which is immediate from
	\begin{align}\label{eq:energyspace}
		\|Ax\|_{\calY_1}^2 = \|Ax\|_{\calY}^2 = \langle Ax,Ax \rangle_{\calY} = h(x,x)=\|x\|_{\calX_h}^2
	\end{align}
	for all $x \in \calX_h$. This gives rise to the notion of \emph{abstract potential energy} and corresponding \emph{energetic space}
	\begin{align*}
		H_p(x) \coloneqq  \tfrac12 h(x,x) = \tfrac12 \|Ax\|^2_{\calY} = \tfrac12 \|x\|^2_{\calX_h},\qquad (\calX_h, \|x\|_{\calX_h}),
	\end{align*}
	see also \cite{PhilReis23}. For the wave equation, this potential energy is given by the $L^2$-norm of the strain variable due to $A=\nabla$.

\begin{figure}
		\centering
		\begin{tikzpicture}
			\draw[thick,black!30!red](7.5,0) circle [radius=2];
			\draw[thick,black!30!green](0.5, -0.5) circle [radius=1];
			\draw[thick,black!30!green](0,0) circle [radius=2];
			\draw[thick,black!30!violet](0.25,-0.25) circle [radius=1.5];
			\node[black!30!green] at (-0.05,-0.25) {$\calX_h$};
			\node[black!30!green] at (-1.25,1) {$\calX_h^*$};
			\node[black] at (-0.75,0.45) {$\calX$};
			\node[black!30!red] at (7.5,0) {$\calY_1$};
			\node[black] at (3.5125,1.8) {$A'$};
			\node[black] at (3.5125,-1.8) {$A$};
			\draw[->, thick, black!50!black, bend left=-25] (6.0,1) to (1.25,1);
			\draw[->, thick, black!50!black, bend left=-25] (1.25,-1) to (6,-1);
		\end{tikzpicture}
		\caption{Illustration of Gelfand triple setting}
		\label{fig:hamiltonian-gelfand-triple}
	\end{figure}
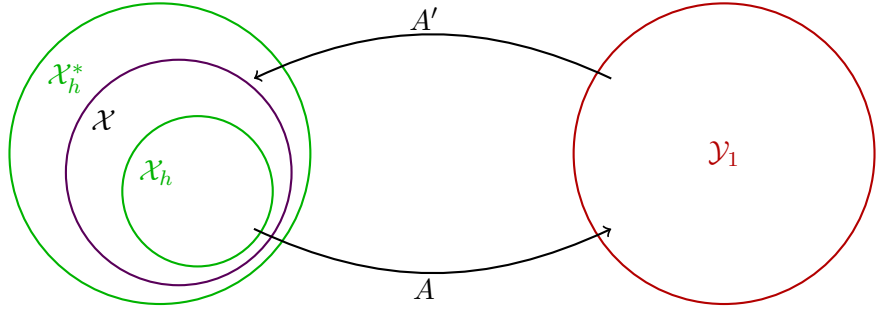

	Note that in \eqref{eq:energyspace}, we just changed the domain and codomain of the operator $A$ and hence, its graph remains the same. In particular, this yields that we are working in a \textit{Gelfand triplet} setting, see Figure \ref{fig:hamiltonian-gelfand-triple}. Observe that by our considerations, the triplet
	\begin{align*}
		\calX_h \overset{A}{\hookrightarrow} \calY_1 \overset{A'}{\hookrightarrow} \calX_h^*
	\end{align*}
	forms a Gelfand triplet, where $A' : \calY_1 \to \calX_h^*$ denotes the Banach space adjoint of $A$ defined in \eqref{eq:bounded_restriction}, in which we identify the intermediate space $\calY \cong \calY_1^*$ with its dual by means of the Riesz isomorphism. The main distinction lies in the fact that, although the triple $(\calX_h , \calX ,\calX_h^*)$ also forms a Gelfand triplet, we employ the isometric isomorphisms $A$ and $A'$ to represent the Riesz isomorphism $S=\calR_{\calX_h} = A'A$ which is the unique bounded extension of the representing operator $S=A^*A$ at the same time. One may think of this as a factorization of the Riesz isomorphism in $\calX_h$ which differs in the choice of the intermediate space.

	\subsection{Boundary triplet for the main operator}\label{subsec:bdd_triplet}
	We now are in the position to define a suitable extension of $\calA$ by extending $A^*$. Then, we define an associated boundary triplet which enables us to later characterize all suitable boundary conditions that render the system scattering or impedance passive. In what follows, we set
	\begin{align*}
		D = 0,\qquad M= I
	\end{align*}
	such that $\calX_M= \calX$ and $\calA \colon \calZ_h \to \calZ_h$ from \eqref{eq:As} simplifies to
	\begin{align}\label{eq:As_neu}
		\dom \calA = \setdef{ \spvek{z_1}{z_2} \in \stack{\calX_h }{\calX_h}}{ S z_1 \in \calX}, \qquad \calA z= \sbmat{
			0}{I}{-S}{0}
		z.
	\end{align}
	The skew-adjointness of $\calA$ is well-known, see \cite{tucsnak2003get}, however we provide the proof for completeness in Lemma~\ref{lem:As_skew} of Appendix~\ref{app:supplementary}.
	We consider the more general case of damping operators $D\in L(\calX)$ with $-D$ dissipative and uniformly positive mass densities $M\in L(\calX)$ in Subsection~\ref{subsec:Mda}.

	To introduce the maximal operator associated to the second order system \eqref{eq:Cauchy_second_order} we assume that there is an extension of $A^*: \calY \supset \dom A^* \to \calX$. Precisely, assume that there is a closed and densely defined linear operator $B: \calY \supset \dom B \to \calX$ such that $(A^*, -B^*)$ builds a dual pair. In other words, \begin{align}\label{eq:dual_pair_ass}
		A^* \subset -B
	\end{align}
	or equivalently, $B^* \subset -A$. Clearly, the product $-BA$ is in general not anymore self-adjoint but it restricts to the self-adjoint operator $S=A^*A$.  In the terminology of the introduction, $S_{\max}=-BA$ is the maximal operator, whose domain carries no boundary condition.
	\begin{rem}\label{rem:choiceA}
	Since the energy space $\calX_h=\dom A$ is fixed, the boundary triplet constructed below can only encode boundary conditions which are not already built into $\calX_h$. For instance, for $A=\grad$ with $\dom A=H_0^1(\Omega)$ (Dirichlet Laplacian), one has $A^*=-\operatorname{div}$ with $\dom A^*=H^{\operatorname{div}}(\Omega)$. Hence, the differential operator $-A^*$ admits no proper extension $B$ of divergence type, the boundary space is trivial, and $\calA_\mathrm{max}=\calA$ is already skew-adjoint. This reflects that the Dirichlet condition is encoded in the energy space $H_0^1(\Omega)$ itself. In order to obtain the maximal freedom of boundary conditions, one therefore chooses $A$ with maximal domain (e.g., $\dom A=H^1(\Omega)$), so that the energy space carries no boundary condition, and all boundary conditions are imposed through the input.
	\end{rem}
	We now aim to identify the operator $B$ as an operator on $\calY_1 = \ran A$. Hence, we define
	\begin{align}\label{eq:B_res}
		B_{\calY_1} : \calY_1 \supset \dom B \cap \calY_1 \to \calX, \qquad  B_{\calY_1}x = Bx
	\end{align}
	for all $x \in \dom{B_{\calY_1}} = \dom B \cap \calY_1$.
	\begin{lem}
		The operator $B_{\calY_1}$ is closed and densely defined.
	\end{lem}
	\begin{proof}
		We first show that $\dom B \cap \calY_1$ is dense in $\calY_1$. Note that taking the closure of a subset of $\calY_1$ is the same as taking the closure with respect to the topology in $\calY$. Now observe that for a closed subset $K \subseteq \calY$ with $K \supseteq \dom B $ and $K\supseteq \calY_1$ we also have
		$K\supseteq \dom B \cap \calY_1$. Hence,
		\begin{align*}
			\overline{\dom B}^{\calY} \cap \overline{\calY_1}^{\calY} &= \left( \bigcap \left\lbrace K \supseteq \dom B \, \middle| \, K \subseteq \calY \ \text{closed} \right\rbrace \right) \bigcap  \left( \bigcap \left\lbrace K \supseteq \ran A \, \middle| \, K \subseteq \calY \ \text{closed} \right\rbrace\right)  \\
			& \subseteq \bigcap \left\lbrace K \supseteq \dom B \cap \calY_1\, \middle| \, K \subseteq \calY \ \text{closed} \right\rbrace  = \overline{\dom B \cap \calY_1}^{\calY}.
		\end{align*}
		But since $B$ is densely defined we have
		\begin{align*}
			\calY_1 =  \overline{\dom B}^{\calY} \cap \overline{\calY_1}^{\calY} \subseteq \overline{\dom B \cap \calY_1}^{\calY} = \overline{\dom B \cap \calY_1}^{\calY_1} \subseteq \calY_1,
		\end{align*}
		hence $\overline{\dom B_{\calY_1}} = \overline{\dom B \cap \calY_1}^{\calY_1} = \calY_1$ such that $B_{\calY_1}$ is densely defined.

		To see the closedness of $B_{\calY_1}$ let $(y_n)_{n \in \N} \subset \dom B_{\calY_1}$ be a sequence with
		\begin{align*}
			y_n \to y \ \text{in}\  \calY_1 , \qquad B_{\calY_1}y_n \to x \ \text{in}\   \calX
		\end{align*}
		as $n \to \infty$. Closedness of $\calY_1$ in $\calY$ and closedness of the operator $B$ yields that $y \in \dom B \cap \calY_1= \dom B_{\calY_1}$ and $B_{\calY_1}y = By=x$.
	\end{proof}
	We now introduce the extension of $\calA$ as introduced in \eqref{eq:As_neu} by the maximal operator $\calA_\mathrm{max}: \calZ_h \to \calZ_h = \stack{\calX_h}{\calX}$ with $\calX_h$ defined in \eqref{eq:defX_h} and where
	\begin{align}\label{eq:defndomA}
		\dom \calA_\mathrm{max}= \left\lbrace \stack{z_1}{z_2} \in\stack{\calX_h}{\calX_h} \, \middle| \, Az_1 \in \dom B_{\calY_1} \right\rbrace , \qquad \calA_\mathrm{max}z= \begin{bmatrix}
			0 & I \\ B_{\calY_1}A & 0
		\end{bmatrix}
		z.
	\end{align}
	At this point, we mention that $\dom \calA_\mathrm{max}= \left\lbrace \stack{z_1}{z_2} \in\stack{\calX_h}{\calX_h}  \, \middle| \, Az_1 \in \dom B \right\rbrace$ since for $z_1 \in \calX_h$ we have that $Az_1 \in \dom B$ if and only if $Az_1 \in  \dom B_{\calY_1} = \dom B \cap \calY_1$.

	In view of the well-known skew-adjointness of $\calA$ (for completeness proven in Lemma~\ref{lem:As_skew} of Appendix~\ref{app:supplementary}) we now collect some properties of our extension $\calA_\mathrm{max}$.
	\begin{lem}\label{lem:A_dense}
		The operator $\calA_\mathrm{max}$ is densely defined, closed and $\calA \subset \calA_\mathrm{max}$. In particular, $\calA_\mathrm{max}^* \subset - \calA_\mathrm{max}$.
	\end{lem}
	\begin{proof}
		We first show that $\calA_\mathrm{max}$ is closed. Let $(z_n)_{n \in \N} \subset \dom \calA_\mathrm{max}$ such that
		\begin{align*}
			(z_n^1,z_n^2) = z_n \to z = (z_1,z_2), \qquad \calA_\mathrm{max}z_n \to y =(y_1, y_2)
		\end{align*}
		in $\calZ_h$ as $n \to \infty$. Obviously, $z_2=y_1 \in \calX_h$ and $z_1 \in \calX_h$ satisfies $Az_1 \in \dom B_{\calY_1}$ by the fact that $B_{\calY_1}Az_n^1 \to y_2$ in $\calX$ combined with the closedness of $B_{\calY_1}$. This shows the closedness of $\calA_\mathrm{max}$. We proceed with the inclusion, $\calA \subset \calA_\mathrm{max}$ which then immediately yields that $\calA_\mathrm{max}$ is densely defined. To this end, let $z \in \dom \calA$. It remains to verify that $Sz_1= A' A z_1 \in \calX$ entails $Az_1 \in \dom B$, where $S:\calX_h \to \calX_h^*$ is to be seen as the Riesz isomorphism, see the discussion after \eqref{eq:defX_h}. Since $\calR_{\calX_h}z_1  = A' A z_1 \in \calX$ we have
		\begin{align*}
			\langle \calR_{\calX_h}z_1 , y \rangle_{\calX} = \langle \calR_{\calX_h}z_1 , y \rangle_{\calX_h^*, \calX_h} = \langle z_1 , y \rangle_{\calX_h} = \langle Az_1 , Ay \rangle_{\calY}
		\end{align*}
		for all $y \in \calX_h=\dom A$. This yields $Az_1 \in \dom A^*$ by definition of the domain of the adjoint operator. As a consequence, $Az_1 \in \dom B$ by \eqref{eq:dual_pair_ass} and since $Az_1 \in \ran A=\calY_1$ is trivially satisfied, we get $Az_1 \in \dom  B_{\calY_1}$. Following the same argumentation and the fact that $\left. B\right|_{\dom A^*} =-A^*$ due to \eqref{eq:dual_pair_ass} we obtain for all $(z_1,z_2)\in \dom \calA$,
		\begin{align*}
			B_{\calY_1} Az_1= BAz_1 = -A^*Az_1
		\end{align*}
		and hence, the actions of $\calA$ and $\calA_\mathrm{max}$ coincide on $\dom \calA$. Eventually, using Lemma~\ref{lem:As_skew},
		\begin{align*}
			\calA_\mathrm{max}^* \subset \calA^* = - \calA \subset -\calA_\mathrm{max}.
		\end{align*}
	\end{proof}
	The preceding Lemma~\ref{lem:A_dense} yields that $\calA_0 \coloneqq  \calA_\mathrm{max}^*$ is closed, densely defined and skew-symmetric. Hence, the prerequisites for a boundary triplet as defined in Definition~\ref{def:operator_bdd_triplet} are met. In this context, we refer to \cite{aigner2025well} where a similar approach, i.e., considering the adjoint operator to obtain a closed, densely defined and skew-symmetric operator to enable the theory of boundary triplets, is chosen for the wave equation.

	In the sequel, we will derive a boundary triplet for the dual pair $(A^*, -B^*)$. Recall that, this means that there exist Hilbert spaces $\calG_1 , \calG_2$ and linear operators
	\begin{align}
		\Lambda = \stack{\Lambda_1}{\Lambda_2}: \dom A \to \stack{\calG_1}{\calG_2^*}, \quad \Pi = \stack{\Pi_1}{\Pi_2}: \dom B \to \stack{\calG_1^*}{\calG_2}
	\end{align}
	such that $(\stack{\calG_1}{\calG_2}, \Lambda,  \Pi)$ satisfies
	\begin{enumerate}
		\item[i)] $\ran \Lambda =\stack{\calG_1}{\calG_2^*}$, $\ran \Pi = \stack{\calG_1^*}{\calG_2}$;
		\item[ii)] the abstract Green identity \begin{align}
			-\langle By,x \rangle_\calX - \langle y, Ax \rangle_\calY = \langle \Pi_1 y, \Lambda_1 x \rangle_{\calG_1^*, \calG_1} - \langle \Pi_2 y, \Lambda_2 x \rangle_{\calG_2, \calG_2^*}
		\end{align}
		holds for all $y \in \dom B, x \in \dom A$.
	\end{enumerate}
	We will see that this assumption ensures that there exists a boundary triplet for $\calA_0^*$ in the sense of Definition \ref{def:operator_bdd_triplet}.
	\begin{thm}\label{thm:bdd_triplet_for_sec_order}
		Let $(\stack{\calG_1}{\calG_2}, \Lambda,  \Pi)$ be a boundary triplet for the dual pair $(A^*, -B^*)$. Define the map $\Gamma=  \stack{\Gamma_0}{\Gamma_1} : \dom \calA_\mathrm{max}= \left\lbrace \stack{z_1}{z_2} \in\stack{\calX_h}{\calX_h} \, \middle| \, Az_1 \in \dom B_{\calY_1} \right\rbrace \to \stack{\calG_1 }{\calG_2} \times \stack{\calG_1 }{\calG_2}^*$ by
		\begin{align*}
			\Gamma_0 z = \begin{bmatrix}
				0 & \Lambda_1 \\ \Pi_2A & 0
			\end{bmatrix} z, \qquad   \Gamma_1 z = \begin{bmatrix}
				-\Pi_1A & 0 \\ 0 & \Lambda_2
			\end{bmatrix} z.
		\end{align*}
		Then, $(\stack{\calG_1 }{\calG_2}, \Gamma_0, \Gamma_1)$ is a boundary triplet for $\calA_0^*=\calA_\mathrm{max}$.
	\end{thm}
	\begin{proof}
		We verify the properties of Definition~\ref{def:operator_bdd_triplet}. First, note that the surjectivity of $\Gamma$ is a direct consequence of the surjectivity of the maps $\Lambda,\Pi$ and $A$.
		Hence, we remain to show the abstract Green identity~\eqref{eq:operator_Green}. For $\roundstack{v}{w}, \roundstack{x}{y} \in \dom \calA_\mathrm{max}$, this follows from
		\begin{align*}
			&\langle \calA_\mathrm{max}\roundstack{v}{w}, \roundstack{x}{y} \rangle_{\calZ_h} + \langle \roundstack{v}{w},  \calA_\mathrm{max} \roundstack{x}{y} \rangle_{\calZ_h} \\
			&= \langle \roundstack{w}{BAv}, \roundstack{x}{y} \rangle_{\calZ_h} + \langle \roundstack{v}{w},   \roundstack{y}{BAx} \rangle_{\calZ_h} \\
			&= \langle Aw, Ax \rangle_\calY  + \langle w, BAx \rangle_\calX + \langle Av, Ay \rangle_\calY + \langle BAv, y \rangle_\calX  \\
			&=\langle \Pi_2 Ax, \Lambda_2 w \rangle_{\calG_2, \calG_2^*} - \langle \Pi_1 Ax, \Lambda_1 w \rangle_{\calG_1^*, \calG_1}   + \langle \Pi_2 Av, \Lambda_2 y \rangle_{\calG_2, \calG_2^*} - \langle \Pi_1 Av, \Lambda_1 y \rangle_{\calG_1^*, \calG_1} \\
			&= \left\langle \Gamma_0\roundstack{v}{w}, \Gamma_1 \roundstack{x}{y}\right\rangle_{ \stack{\calG_1 }{\calG_2}, \stack{\calG_1 }{\calG_2}^*}+ \left\langle \Gamma_1 \roundstack{v}{w}, \Gamma_0 \roundstack{x}{y}\right\rangle_{ \stack{\calG_1 }{\calG_2}^* ,  \stack{\calG_1 }{\calG_2}}.
		\end{align*}
	\end{proof}

	\subsection{Abstract Helmholtz decomposition and equivalence transforms}\label{subsec:Helmholtz}

	In this subsection, we transform the operator $\calA_\mathrm{max}$ using the isometric isomorphism $A: \calX_h \to \calY_1$. The motivation for this is to provide an abstract operator-theoretic framework to transformations between different formulations of second-order equations in the literature, in particular in view of port-Hamiltonian systems, see the wave equation of \eqref{eq:wavebeispiel}.

	First, recall the splitting of the space $\calY$ into orthogonal subspaces
	\begin{align}\label{eq:Helmholtz}
		\calY = \calY_0 \oplus \calY_1 = (\ran A )^\perp \oplus \overline{\ran A} = \ker A^* \oplus \ran A,
	\end{align}
	which may be seen as an \textit{abstract Helmholtz decomposition}.
	This nomenclature is motivated by the example $A = \nabla$, where the above sum models the decomposition of $L^2$ into gradient fields and divergence-free vector fields.

	Define the equivalence transformation $\calB_1: \stack{\calY_1}{\calX}  \supset \dom \calB_1 \to \stack{\calY_1}{\calX}$ where
	\begin{align}\label{eq:defndomB_1}
		\calB_1 = \begin{bmatrix}
			A & 0 \\ 0 & I
		\end{bmatrix} \calA_\mathrm{max}\begin{bmatrix}
			A^{-1} & 0 \\ 0 & I
		\end{bmatrix} = \begin{bmatrix}
			0 & A \\ B_{\calY_1} & 0
		\end{bmatrix}, \qquad \dom \calB_1 = \begin{bmatrix}
			A & 0 \\ 0 & I
		\end{bmatrix} \dom \calA_\mathrm{max}= \stack{ \dom B_{\calY_1} }{ \calX_h}.
	\end{align}
	The following lemma makes explicit the relationship between the boundary triplet for the original system $\calA_0^*$ and its transformed variant $\calB_1$ .
	It is immediate from the definition of the equivalence transforms that $\calB_1$ is closed, densely defined and the adjoint of a skew-symmetric operator.
	\begin{lem}\label{lem:bdd_triplet_rel}
		Let $(\stack{\calG_1 }{\calG_2}, \Gamma_0, \Gamma_1)$ be the boundary triplet for $\calA_0^*=\calA_\mathrm{max}$ defined in Theorem \ref{thm:bdd_triplet_for_sec_order}.
		Define the map $\Xi= \stack{\Xi_0}{\Xi_1} : \dom \calB_1 \to \stack{\calG_1 }{\calG_2} \times \stack{\calG_1 }{\calG_2}^*$ via
		\begin{align*}
			\Xi_0 x = \begin{bmatrix}
				0 & \Lambda_1 \\ \Pi_2 & 0
			\end{bmatrix} x, \qquad   \Xi_1 x = \begin{bmatrix}
				-\Pi_1 & 0 \\ 0 & \Lambda_2
			\end{bmatrix} x.
		\end{align*}
		Then, $(\stack{\calG_1 }{\calG_2}, \Xi_0, \Xi_1)$ is a boundary triplet for $\calB_1$. In particular, $\Xi_i = \Gamma_i \left[ \begin{smallmatrix}
			A^{-1} & 0 \\ 0 & I
		\end{smallmatrix}\right]$ for $i=1,2$.
	\end{lem}
	\begin{proof}
		As surjectivity of the boundary maps is immediate, it is enough to check is the validity of the abstract Green identity:
		\begin{align*}
			&\langle \calB_1 \roundstack{v}{w}, \roundstack{x}{y} \rangle_{\stack{\calY_1}{\calX}} + \langle \roundstack{v}{w},  \calB_1  \roundstack{x}{y} \rangle_{\stack{\calY_1}{\calX}} \\
			&=\langle \left[ \begin{smallmatrix}
				A & 0 \\ 0 & I
			\end{smallmatrix}\right] \calA_0^* \left[ \begin{smallmatrix}
				A^{-1} & 0 \\ 0 & I
			\end{smallmatrix}\right] \roundstack{v}{w}, \roundstack{x}{y} \rangle_{\stack{\calY_1}{\calX}} + \langle \roundstack{v}{w},   \left[ \begin{smallmatrix}
				A & 0 \\ 0 & I
			\end{smallmatrix}\right] \calA_0^* \left[ \begin{smallmatrix}
				A^{-1} & 0 \\ 0 & I
			\end{smallmatrix}\right]  \roundstack{x}{y} \rangle_{\stack{\calY_1}{\calX}} \\
			&=\langle \calA_0^* \left[ \begin{smallmatrix}
				A^{-1} & 0 \\ 0 & I
			\end{smallmatrix}\right] \roundstack{v}{w}, \left[ \begin{smallmatrix}
				A' & 0 \\ 0 & I
			\end{smallmatrix}\right] \roundstack{x}{y} \rangle_{\calZ_h, \calZ_h^*} + \langle  \left[ \begin{smallmatrix}
				A' & 0 \\ 0 & I
			\end{smallmatrix}\right] \roundstack{v}{w},  \calA_0^* \left[ \begin{smallmatrix}
				A^{-1} & 0 \\ 0 & I
			\end{smallmatrix}\right]  \roundstack{x}{y} \rangle_{\calZ_h^*, \calZ_h} \\
			&=\langle \calA_0^* \left[ \begin{smallmatrix}
				A^{-1} & 0 \\ 0 & I
			\end{smallmatrix}\right] \roundstack{v}{w}, \left[ \begin{smallmatrix}
				\calR_{\calX_h} A^{-1} & 0 \\ 0 & I
			\end{smallmatrix}\right] \roundstack{x}{y} \rangle_{\calZ_h, \calZ_h^*} + \langle  \left[ \begin{smallmatrix}
				\calR_{\calX_h} A^{-1} & 0 \\ 0 & I
			\end{smallmatrix}\right] \roundstack{v}{w},  \calA_0^* \left[ \begin{smallmatrix}
				A^{-1} & 0 \\ 0 & I
			\end{smallmatrix}\right]  \roundstack{x}{y} \rangle_{\calZ_h^*, \calZ_h} \\
			&=\langle \calA_0^* \left[ \begin{smallmatrix}
				A^{-1} & 0 \\ 0 & I
			\end{smallmatrix}\right] \roundstack{v}{w}, \left[ \begin{smallmatrix}
				A^{-1} & 0 \\ 0 & I
			\end{smallmatrix}\right] \roundstack{x}{y} \rangle_{\calZ_h} + \langle  \left[ \begin{smallmatrix}
				A^{-1} & 0 \\ 0 & I
			\end{smallmatrix}\right] \roundstack{v}{w},  \calA_0^* \left[ \begin{smallmatrix}
				A^{-1} & 0 \\ 0 & I
			\end{smallmatrix}\right]  \roundstack{x}{y} \rangle_{ \calZ_h} \\
			&= \left\langle \Gamma_0 \left[ \begin{smallmatrix}
				A^{-1} & 0 \\ 0 & I
			\end{smallmatrix}\right] \roundstack{v}{w}, \Gamma_1  \left[ \begin{smallmatrix}
				A^{-1} & 0 \\ 0 & I
			\end{smallmatrix}\right] \roundstack{x}{y}\right\rangle_{\stack{\calG_1 }{\calG_2}, \stack{\calG_1 }{\calG_2}^*}+ \left\langle \Gamma_1 \left[ \begin{smallmatrix}
				A^{-1} & 0 \\ 0 & I
			\end{smallmatrix}\right] \roundstack{v}{w}, \Gamma_0 \left[ \begin{smallmatrix}
				A^{-1} & 0 \\ 0 & I
			\end{smallmatrix}\right] \roundstack{x}{y}\right\rangle_{\stack{\calG_1 }{\calG_2}^*, \stack{\calG_1 }{\calG_2}} \\
			&= \left\langle \Xi_0 \roundstack{v}{w}, \Xi_1 \roundstack{x}{y}\right\rangle_{\stack{\calG_1 }{\calG_2}, \stack{\calG_1 }{\calG_2}^*}+ \left\langle \Xi_1 \roundstack{v}{w}, \Xi_0 \roundstack{x}{y}\right\rangle_{\stack{\calG_1 }{\calG_2}^*, \stack{\calG_1 }{\calG_2}}
		\end{align*}
		for all $\roundstack{v}{w}, \roundstack{x}{y} \in \dom \calB_1$. Hence, the assertion is proved.
	\end{proof}
	In many applications, the operator $\calB_1$ is seen as an operator in $\stack{\calY}{\calX}$ instead of just $\stack{\calY_1}{\calX} = \stack{\ran A}{\calX}$. This is possible due to the abstract Helmholtz decomposition~\eqref{eq:Helmholtz} of $\calY$ and the symmetry assumption $B^* \subset -A$. More precisely, define the operator
	\begin{align}\label{eq:ende}
		\calB  :\stack{\calY}{\calX} \supset \dom \calB \to \stack{\calY}{\calX} , \qquad \dom \calB = \stack{ \dom B}{ \dom A}, \qquad  \calB = \begin{bmatrix}
			0 & A \\ B & 0
		\end{bmatrix}
	\end{align}
	which is the extension of $\calB_1$ to the whole space $\stack{\calY}{\calX}$. This follows from
	the fact that $\ker A^* \subset \ker B$ such that going from $B_{\calY_1}$ to $B$ in the second row of $\calB$ does not affect the action of the operator. In particular, the graphs of $\calB_1$ and $\calB$ in $\stack{\calY}{\calX} \times \stack{\calY}{\calX}$ coincide.

	\subsection{The associated boundary node}\label{subsec:bdd_node}
	We now define the boundary node corresponding to $\calA_\mathrm{max}$ defined in \eqref{eq:defndomA}. The boundary triplet $(\calG, \Gamma_0, \Gamma_1)$ for $\calA_0^*$ derived in Theorem~\ref{thm:bdd_triplet_for_sec_order} allows us to define a wide range of suitable input and output boundary operators to define a boundary control system of the form \eqref{eq:ODEnode}, in particular in view of scattering or impedance passivity. At first glance a suitable choice would be $G=\Gamma_0$ and $K=\Gamma_1$ (or vice versa) as performed in~\cite{arov2012boundary}, which in case of the wave equation corresponds to Dirichlet and Neumann boundary control and observation. Precisely,
	\begin{align}\label{eq:natural_bc}
		\Theta = \left(\begin{bsmallmatrix}
			\Gamma_0 \\ \calA_0^* \\ \Gamma_1
		\end{bsmallmatrix}, \begin{bsmallmatrix}
			\calG \vphantom{\Gamma_0} \\ \calZ_h \vphantom{\calA_0^*}  \\ \calG^* \vphantom{\Gamma_1}
		\end{bsmallmatrix}\right)
	\end{align}
	with domain $\dom \Theta=\dom \calA_0^*$ is an internally well-posed boundary node, see \cite[Thm 4.8]{arov2012boundary}. Its main operator even generates a unitary semigroup since $\calA_0^*$ restricted to $\ker \Gamma_0$ is skew-adjoint. In this paper, we generalize this result by allowing Robin-type boundary controls
	\begin{align}\label{eq:Robin_control}
		u(t)= \tfrac{1}{2}(1-P)\calR_{\calG} \Gamma_0z(t) + \tfrac{1}{2} (1+P)\Gamma_1z(t) ,
	\end{align}
	where $P$ is some contraction on $\calG^*$. We will observe that such input maps arise when the external Cayley transforms of specific scattering passive boundary nodes are considered. We also observe that Robin-type input maps as in \eqref{eq:Robin_control} exploit the full freedom of choice of boundary conditions that lead to maximal dissipative main operators.
	\begin{thm}\label{thm:main}
		Let $(\calG, \Gamma_0, \Gamma_1)$ be the boundary triplet for $\calA_0^*$ provided in Theorem \ref{thm:bdd_triplet_for_sec_order}. Let $P \in L(\calG^*)$ and
		define the colligation
		\begin{align*}
			N_{\calZ_h} = \left(\begin{bsmallmatrix}
				\tfrac{1}{\sqrt{2}}(\calR_{\calG}\Gamma_0 + \Gamma_1) \\ \calA_0^* \\ \tfrac{-1}{\sqrt{2}}P(\calR_{\calG}\Gamma_0 - \Gamma_1)
			\end{bsmallmatrix}, \begin{bsmallmatrix}
				\calG^* \vphantom{\tfrac{1}{\sqrt{2}}} \vphantom{\Gamma_0} \\ \calZ_h \vphantom{\calA_0^*}  \\ \calG^* \vphantom{\Gamma_1} \vphantom{\tfrac{1}{\sqrt{2}}}
			\end{bsmallmatrix}\right)
		\end{align*}
		with domain $\dom N_{\calZ_h}=\dom \calA_0^*$. Then, $N_{\calZ_h}$ defines a scattering passive boundary node if and only if $\|P\|_{L(\calG^*)}\leq 1$. Moreover, $N_{\calZ_h}$ defines a scattering energy-preserving boundary node if and only if $P$ is isometric, i.e., $\|Pg\|_{\calG^*} = \|g\|_{\calG^*}$ for all $g \in \calG^*$.
	\end{thm}
	\begin{proof}
		First, observe that clearly $N_{\calZ_h}$ is a colligation.
		We show that $\|P\|_{L(\calG^*)}\leq 1$ implies scattering passivity of $N_{\calZ_h}$ via \cite[Theorem 2.5]{malinen2007impedance}.
		To this end, we have to verify that
		\begin{itemize}
			\item[(i)] $\stack{\tfrac{1}{\sqrt{2}}(\calR_{\calG}\Gamma_0 + \Gamma_1)}{\alpha - \calA_0^*}$ is surjective for all $\alpha \in \C_+$;
			\item[(ii)] for all $z \in \dom \calA_0^*$ we have
			\begin{align}\label{eq:pass_ineq_A_0_star}
				2\Re \langle z, \calA_0^*z \rangle_\calX + \frac{1}{2}||P(\calR_{\calG}\Gamma_0 - \Gamma_1)z||^2_{\calG^*} \leq \frac{1}{2}||(\calR_{\calG}\Gamma_0 + \Gamma_1)z||^2_{\calG^*}.
			\end{align}
		\end{itemize}
		In view of the discussion after \cite[Theorem 2.5]{malinen2007impedance}, $(i)$ holds if and only if $\tfrac{1}{\sqrt{2}}(\calR_{\calG}\Gamma_0 + \Gamma_1)$ and $\alpha - A_{\calZ_h}$ are surjective where $A_{\calZ_h}=\left. \calA_0^* \right|_{\ker \calR_{\calG}\Gamma_0 + \Gamma_1}$.
		To see this, first note that clearly $\nextto{\calR_{\calG}}{1}: \stack{\calG}{\calG^*} \to \calG^*$ is surjective.
		Consequently, using also surjectivity of $\Gamma$, for all $u \in \calG^*$, there exists some $z \in \dom \calA_0^*$ with $\nextto{\calR_{\calG}}{1}\Gamma z=u$.
		Furthermore, the operator $A_{\calZ_h}$ generates a contraction semigroup on $\calZ_h$ by Theorem \ref{thm:former_Kurula_Zwart} and the choice $P=0$. Then, the Lumer-Philipps Theorem \cite[Thm. 3.8.4]{tucsnak2009observation} yields that $A_{\calZ_h}$ is maximal dissipative and hence $\C_+ \subset \rho( A_{\calZ_h})$ which proves in particular surjectivity of $\alpha - A_{\calZ_h}$ for all $\alpha \in \C_+$. \\ It remains to prove $(ii)$.
		From the polarization identity in Hilbert spaces, and as $P$ is a contraction, we obtain for all $z \in \dom \calA_0^*$
		\begin{align*}
			2 \Re \langle \calR_\calG \Gamma_0 z, \Gamma_1 z \rangle_{\calG^*}& =  \frac{1}{2}||(\calR_{\calG}\Gamma_0 + \Gamma_1)z||^2_{\calG^*} - \frac{1}{2}||(\calR_{\calG}\Gamma_0 - \Gamma_1)z||^2_{\calG^*} \\
			& \leq \frac{1}{2}||(\calR_{\calG}\Gamma_0 + \Gamma_1)z||^2_{\calG^*} - \frac{1}{2}||P(\calR_{\calG}\Gamma_0 - \Gamma_1)z||^2_{\calG^*}.
		\end{align*}
		As $(\calG, \Gamma_0, \Gamma_1)$ is a boundary triplet for $\calA_0^*$, the abstract Green identity \eqref{eq:operator_Green} is satisfied. As a consequence,
		\begin{align}\label{eq:inequality_scatt}
			\begin{split}
				2\Re \langle z, \calA_0^*z \rangle_\calX + \frac{1}{2}||P(\calR_{\calG}\Gamma_0 - \Gamma_1)z||^2_{\calG^*} & =  2 \Re \langle \calR_\calG \Gamma_0 z, \Gamma_1 z \rangle_{\calG^*}  +  \frac{1}{2}||P(\calR_{\calG}\Gamma_0 - \Gamma_1)z||^2_{\calG^*} \\
				& \leq  \frac{1}{2}||(\calR_{\calG}\Gamma_0 + \Gamma_1)z||^2_{\calG^*}
			\end{split}
		\end{align}
		for all $z \in \dom \calA_0^*$. Hence, $N_{\calZ_h}$ is a scattering passive boundary node. In case $P$ is an isometric operator on the boundary space $\calG^*$, then the inequality in \eqref{eq:inequality_scatt} becomes an equality such that $N_{\calZ_h}$ is scattering energy-preserving.

		Conversely, let $N_{\calZ_h}$ be a scattering passive boundary node which in particular yields \eqref{eq:pass_ineq_A_0_star}. Since $z \mapsto (\calR_{\calG}\Gamma_0 - \Gamma_1)z$ is surjective, any $y \in \calG^*$ can be written as $y \coloneqq  (\calR_{\calG}\Gamma_0 - \Gamma_1)z$ for a particular $z \in \dom \calA_0^*$. Then, the passivity inequality \eqref{eq:pass_ineq_A_0_star} yields, for all $y\in \calG^*$,
		\begin{align*}
			\frac{1}{2} || P y||_{\calG^*}^2 & \leq  \frac{1}{2}||y  + 2\Gamma_1z||^2_{\calG^*} -   2\Re \langle z, \calA_0^*z \rangle_\calX \\
			&=  \frac{1}{2}||y ||^2_{\calG^*} + 2 \Re \langle y, \Gamma_1 z \rangle_{\calG^*}   + 2||\Gamma_1z||^2_{\calG^*} - 2 \Re \langle \calR_\calG \Gamma_0 z, \Gamma_1 z \rangle_{\calG^*} \\
			&=  \frac{1}{2}||y ||^2_{\calG^*} + 2 \Re \langle y, \Gamma_1 z \rangle_{\calG^*}   + 2||\Gamma_1z||^2_{\calG^*} - 2 \Re \langle y+\Gamma_1 z, \Gamma_1 z \rangle_{\calG^*}  \\
			&=   \frac{1}{2}||y ||^2_{\calG^*}
		\end{align*}
		such that $P$ is a contraction.
		If $N_{\calZ_h}$ is scattering energy-preserving, the above chain of inequalities holds with equality, such that $P$ is isometric.
	\end{proof}

	\begin{cor}
		Let $N_{\calZ_h}$ be defined as in Theorem \ref{thm:main} and consider its external Cayley transform with parameter $\beta = 1$
		\begin{align}\label{eq:ext_Cayley_Z_h}
			N_{\calZ_h}^{(1)} = \left(\begin{bsmallmatrix}
				\tfrac{1}{2}((1-P)\calR_{\calG} \Gamma_0+ (1+P)\Gamma_1) \\ \calA_0^* \\ \tfrac{1}{2}((1+P)\calR_{\calG} \Gamma_0+(1-P)\Gamma_1)
			\end{bsmallmatrix}, \begin{bsmallmatrix}
				\calG^* \vphantom{\tfrac{1}{\sqrt{2}}} \vphantom{\Gamma_0} \\ \calZ_h \vphantom{\calA_0^*}  \\ \calG^* \vphantom{\Gamma_1} \vphantom{\tfrac{1}{\sqrt{2}}}
			\end{bsmallmatrix}\right)
		\end{align}
		with domain $\dom N_{\calZ_h}^{(1)}=\dom \calA_0^*$. If $P$ is a contraction, $N_{\calZ_h}^{(1)}$ is an impedance passive colligation and if $P$ is unitary, $N_{\calZ_h}^{(1)}$ is an impedance energy-preserving colligation.
	\end{cor}
	\begin{proof}
		The external Cayley transform trivially satisfies the equality
		\begin{align*}
			\left(N_{\calZ_h}^{(1)}\right)^{(1)} = N_{\calZ_h}.
		\end{align*}
		Then, combining Theorem \ref{thm:main} and Definition \ref{def:imped_bn} yields the result.
	\end{proof}
	Note that the rather simple choice of Neumann and Dirichlet-like input and output boundary maps as in \eqref{eq:natural_bc} can be recovered by setting $P=\pm 1$ in the preceding Corollary.
	\begin{rem}
		The authors emphasize that the input boundary operator of $N_{\calZ_h}$ is natural and maximal in the following sense: the main operator of an impedance passive internally well-posed boundary node is always the generator of contraction semigroup on the state space, cf.\ \cite[Proposition 4.1]{malinen2007impedance}. In our setting, we derive a boundary triplet $(\calG, \Gamma_0, \Gamma_1)$ for the operator $\calA_0^*$ corresponding to the second-order system \eqref{eq:second_order}. Following Theorem \ref{thm:former_Kurula_Zwart}, for any restriction $B$ of $\calA_0^*$ that generates a contractive semigroup there exists a contraction $P \in L(\calG^*)$ such that $\dom B= \ker (P-1)\calR_\calG \Gamma_0 - (P+1)\Gamma_1$ and this is exactly the kernel of the input boundary operator of $N_{\calZ_h}$.
	\end{rem}
	\begin{cor}\label{cor:well_posed}
		Let $N_{\calZ_h}$ be defined as in Theorem \ref{thm:main} with contractive $P \in L(\calG^*)$. Then its external Cayley transform $N_{\calZ_h}^{(1)}$ defined in \eqref{eq:ext_Cayley_Z_h} is an impedance passive internally well-posed boundary node.
	\end{cor}
	\begin{proof}
		Following the argumentation in the beginning of the proof Theorem \ref{thm:main} shows that the conditions of \cite[Proposition 4.1]{malinen2007impedance} are satisfied. This yields the claim.
	\end{proof}
	\subsection{Nontrivial density and dissipation operators}\label{subsec:Mda}
	We now also allow for nonzero dissipation operators and general mass density operators. To this end, consider $M \in L(\calX)$ self-adjoint and uniformly positive, i.e. there exists a constant $\alpha > 0$ such that
	\begin{align*}
		\|Mx\|_\calX \geq \alpha \|x\|_\calX
	\end{align*}
	for all $x \in \calX$. The damping operator $D\in L(\calX)$ is such that $-D$ is dissipative. Recall the definition of the weighted Hilbert space $\calX_M$ as $\calX$ endowed with $\langle M \cdot, \cdot \rangle_\calX$. We slightly modify $\calA_\mathrm{max}$ by introducing the operator $\calA_\mathrm{max,M}: \stack{\calX_h}{\calX_M} \supset \dom \calA_\mathrm{max,M} \to\stack{\calX_h}{\calX_M}$ where
	\begin{align*}
		\dom \calA_\mathrm{max,M} =\left\lbrace \stack{z_1}{z_2} \in\stack{ \calX_h }{ M\calX_h } \, \middle| \, Az_1 \in \dom B_{\calY_1} \right\rbrace  , \qquad \calA_\mathrm{max,M} = \begin{bmatrix}
			0 & M^{-1} \\ B_{\calY_1}A & \ -DM^{-1}
		\end{bmatrix} .
	\end{align*}
	Note that $\calA_\mathrm{max,M}$ is just the product
	\begin{align*}
		\calA_\mathrm{max,M}= \left(\calA_\mathrm{max}- \begin{bmatrix}
			0 & 0 \\ 0 & D
		\end{bmatrix} \right) \begin{bmatrix}
			I & 0 \\ 0 & M^{-1}
		\end{bmatrix}
	\end{align*}
	considered in the weighted Hilbert space $\stack{\calX_h}{\calX_M}$.
	\begin{thm}\label{thm:main2}
		Let $(\calG, \Gamma_0, \Gamma_1)$ be the boundary triplet for $\calA_0^*$ provided in Theorem \ref{thm:bdd_triplet_for_sec_order}. Let $P \in L(\calG^*)$ with $\|P\|_{L(\calG^*)}\leq 1$. Define the colligation
		\begin{align*}
            N= \left(\begin{bsmallmatrix}
	\tfrac{1}{2}((1-P)\calR_{\calG} \Gamma_0+ (1+P)\Gamma_1) \left[ \begin{smallmatrix}
		I & 0 \\ 0 & M^{-1}
	\end{smallmatrix}\right]  \\ \calA_\mathrm{max,M} \vphantom{\stack{\calX_h}{\calX_M}} \\ \tfrac{1}{2}((1+P)\calR_{\calG} \Gamma_0+ (1-P)\Gamma_1) \left[ \begin{smallmatrix}
		I & 0 \\ 0 & M^{-1}
	\end{smallmatrix}\right]
\end{bsmallmatrix}, \begin{bsmallmatrix}
	\calG^* \vphantom{\left[ \begin{smallmatrix} I & 0 \\ 0 & M^{-1} \end{smallmatrix}\right]} \\ \stack{\calX_h}{\calX_M}  \\ \calG^* \vphantom{\left[ \begin{smallmatrix} I & 0 \\ 0 & M^{-1} \end{smallmatrix}\right]}
\end{bsmallmatrix}\right)
		\end{align*}
		with domain $\dom N=\dom \calA_\mathrm{max,M}$. Then, $N$ defines an impedance passive internally well-posed boundary node.
	\end{thm}
	\begin{proof}
		It is shown in \cite[Lem. 2.14]{gernandt2024stability} that if $(\calG, \Gamma_0,  \Gamma_1)$ defines a boundary triplet for $\calA_\mathrm{max}$, then also $(\calG, \Gamma_0 \left[\begin{smallmatrix}
			I & 0 \\ 0 & M^{-1}
		\end{smallmatrix}\right], \Gamma_1 \left[\begin{smallmatrix}
			I & 0 \\ 0 & M^{-1}
		\end{smallmatrix}\right])$ is a boundary triplet for $\calA_\mathrm{max}\left[\begin{smallmatrix}
			I & 0 \\ 0 & M^{-1}
		\end{smallmatrix}\right]$ as an operator in $\stack{\calX_h}{\calX_M}$. Then, it is a direct consequence of Theorem \ref{thm:main} and Corollary \ref{cor:well_posed} that $N$ defines an impedance passive internally well-posed boundary node by the fact that the bounded perturbation $\left[ \begin{smallmatrix}
			0 & 0 \\ 0 & -D
		\end{smallmatrix}\right]$ and the change of the inner product do not affect the properties in Definition \ref{def:boundarynode} and the passivity inequality \eqref{eq:pass_ineq_A_0_star}.
	\end{proof}
	The isometric isomorphism $A:\calX_h \to \calY_1$ gives rise to an equivalence transform of the boundary node $N$. Note that the following theorem is in the spirit of Subsection \ref{subsec:Helmholtz}. To this end, define the similarity transform in $\stack{\calY_1}{\calX_M} $:
	\begin{align*}
		\calB_M \coloneqq  \left[ \begin{smallmatrix}
			A & 0 \\ 0 & I_{\calX_M}
		\end{smallmatrix}\right] \calA_\mathrm{max,M}  \left[ \begin{smallmatrix}
			A^{-1} & 0 \\ 0 & I_{\calX_M}
		\end{smallmatrix}\right] = \begin{bmatrix}
			0 & AM^{-1} \\ B_{\calY_1} & \ -DM^{-1}
		\end{bmatrix} , \qquad \dom \calB_M \coloneqq  \left[ \begin{smallmatrix}
			A & 0 \\ 0 & I_{\calX_M}
		\end{smallmatrix}\right] \dom \calA_\mathrm{max,M}.
	\end{align*}

	Recall the definition of the boundary triplet $(\calG, \Xi_0, \Xi_1)$ for $\calB_1$ as in Lemma \ref{lem:bdd_triplet_rel}.
	\begin{thm}\label{thm:trafo}
		Let $P \in L(\calG^*)$ with $\|P\|_{L(\calG^*)}\leq 1$. Define the colligation
		\begin{equation}\label{eq:trafo1}
			L \coloneqq   \left(\begin{bsmallmatrix}
				\tfrac{1}{2}((1-P)\calR_{\calG} \Xi_0+ (1+P)\Xi_1)\left[ \begin{smallmatrix}
					I_{\calY_1} & 0 \\ 0 & M^{-1}
				\end{smallmatrix}\right]   \\ \calB_M  \vphantom{\stack{\calY_1}{\calX_M} }\\ \tfrac{1}{2}((1+P)\calR_{\calG} \Xi_0+ (1-P)\Xi_1)\left[ \begin{smallmatrix}
					I_{\calY_1} & 0 \\ 0 & M^{-1}
				\end{smallmatrix}\right]
			\end{bsmallmatrix}, \begin{bsmallmatrix}
				\calG^* \vphantom{\left[ \begin{smallmatrix}
					I_{\calY_1} & 0 \\ 0 & M^{-1}
				\end{smallmatrix}\right]} \\ \stack{\calY_1}{\calX_M}  \\ \calG^* \vphantom{\left[ \begin{smallmatrix}
					I_{\calY_1} & 0 \\ 0 & M^{-1}
				\end{smallmatrix}\right]}
			\end{bsmallmatrix}\right) .
		\end{equation}
		with domain $\dom L = \dom \calB_M$. Then, $L$ is an impedance passive internally well-posed boundary node.
	\end{thm}
	\begin{proof}
		The first part of the proof is a direct consequence of the fact that $N$ is already a boundary node and that $L$ is the similarity transform of $N$ with similarity $\left[ \begin{smallmatrix}
			A & 0 \\ 0 & I_{\calX_M}
		\end{smallmatrix}\right]$. It is well-known that similarity transforms of boundary nodes again define boundary nodes, see e.g. \cite[Ex. 2.3.7]{Staffans2005}. The impedance passivity follows using the same argumentation as in the proof of Theorem \ref{thm:main}.
	\end{proof}
	\begin{rem}
		In view of the discussion at then end of Subsection~\ref{subsec:Helmholtz}, see \eqref{eq:ende}, the boundary node $L$ given in Theorem \ref{thm:trafo} may also be extended to a boundary node $L$ on the state space $\stack{\calY}{\calX_M}$.
	\end{rem}
	The ensuing corollary is immediate from Lemma~\ref{lem:solex} and the fact that the impedance passive boundary node $N$ defined in Theorem \ref{thm:main2} is internally well-posed.
	\begin{cor}\label{cor:solex_pH}
		Let $P \in L(\calG^*)$ with $\|P\|_{L(\calG^*)}\leq 1$. Let $z_0 \in \stack{\calX_h}{\calX_M}$ and $u \in C^2(\R_{\geq 0};\calG^*)$ with $\tfrac{1}{2}((1-P)\calR_{\calG} \Gamma_0+ (1+P)\Gamma_1) \left[ \begin{smallmatrix}
			I & 0 \\ 0 & M^{-1}
		\end{smallmatrix}\right] z_0=u(0)$. Then, there exists a unique smooth solution $z \in  C^1(\R_{\geq 0} ; \stack{\calX_h}{\calX_M}) \cap C(\R_{\geq 0} ; \dom \calA_\mathrm{max,M})$ of
		\begin{align}\label{eq:pH_dynamics}
			u(t)&=\tfrac{1}{2}((1-P)\calR_{\calG} \Gamma_0+ (1+P)\Gamma_1) \left[ \begin{smallmatrix}
				I & 0 \\ 0 & M^{-1}
			\end{smallmatrix}\right]z(t), \nonumber \\ \dot{z}(t)& =\begin{bmatrix}
				0 & M^{-1} \\ B_{\calY_1}A & \ -DM^{-1}
			\end{bmatrix}  z(t),  \\
			y(t)&=\tfrac{1}{2}((1+P)\calR_{\calG} \Gamma_0+ (1-P)\Gamma_1) \left[ \begin{smallmatrix}
				I & 0 \\ 0 & M^{-1}
			\end{smallmatrix}\right]z(t) \nonumber
		\end{align} with $z(0)=z_0$.
	\end{cor}
     \begin{rem}\label{rem:L2wellposed}
    	 The regularity $u\in C^2$ and the compatibility condition $Gz_0=u(0)$ in the  preceding corollary are  standard assmuptions for existence of classical solutions. They ensure that the input equation $u(t)=Gz(t)$ holds pointwise in time. In the framework of system nodes even $u\in W^{2,1}$ suffices, see \cite[Thm.~4.3.9]{Staffans2005}. The solutions corresponding to non-smooth inputs are then generalized in the sense that they are limits of classical solutions, see e.g. \cite[Def.~7]{PhilReis23}. For the scattering passive node $N_{\calZ_h}$ of Theorem~\ref{thm:main}, integration of \eqref{eq:passivity_ineq} gives
    	\begin{equation*}
        	\|z(t)\|^2_{\calZ_h}+\|y\|^2_{L^2(0,t;\calG^*)}\le \|z_0\|^2_{\calZ_h}+\|u\|^2_{L^2(0,t;\calG^*)},\qquad t\ge0,
    	\end{equation*}
    	so $N_{\calZ_h}$ is \textit{well-posed} in the sense of \cite[Def.~10]{PhilReis23}. Precisely, every $z_0\in\calZ_h$ and $u\in L^2_{\mathrm{loc}}(\R_{\ge0};\calG^*)$ yield a unique generalized solution. In contrast, an impedance passive system is well-posed if and only if its transfer function is bounded on some vertical line in $\C_+$, cf.\ \cite[Prop.~30]{PhilReis23}.
    \end{rem}
	\begin{cor}
		Let $P \in L(\calG^*)$ with $\|P\|_{L(\calG^*)}\leq 1$. Let $(u,z,y) \in C^2(\R_{\geq 0};\calG^*) \times C^1(\R_{\geq 0} ; \stack{\calX_h}{\calX_M}) \cap C(\R_{\geq 0} ; \dom \calA_\mathrm{max,M}) \times C(\R_{\geq 0};\calG^*)$ be the unique smooth solution of \eqref{eq:pH_dynamics} with $z(0)=z_0$. The unique classical trajectory
		$(u, w, v) \in C^2(\R_{\geq 0};\calG^*) \times C^1(\R_{\geq 0} ; \stack{\calY_1}{\calX_M}) \cap C(\R_{\geq 0} ; \dom \calB_M) \times C(\R_{\geq 0};\calG^*)$ that solves
		\begin{align}\label{eq:pH_dynamics_2}
			u(t)&=\tfrac{1}{2}((1-P)\calR_{\calG} \Xi_0+ (1+P)\Xi_1) \left[ \begin{smallmatrix}
				I_{\calY_1} & 0 \\ 0 & M^{-1}
			\end{smallmatrix}\right] w(t), \nonumber \\ \dot{w}(t)& =\begin{bmatrix}
				0 & AM^{-1} \\ B_{\calY_1} & \ -DM^{-1}
			\end{bmatrix}  w(t),  \\
			v(t)&=\tfrac{1}{2}((1+P)\calR_{\calG} \Xi_0+ (1-P)\Xi_1) \left[ \begin{smallmatrix}
				I_{\calY_1} & 0 \\ 0 & M^{-1}
			\end{smallmatrix}\right] w(t) \nonumber
		\end{align} with $w(0)=\left[ \begin{smallmatrix}
			A & 0 \\ 0 & I
		\end{smallmatrix}\right]z_0$
		Then,
		\[
		(u,w, v)=(u, \left[ \begin{smallmatrix}
			A & 0 \\ 0 & I
		\end{smallmatrix}\right]z, y).
		\]
	\end{cor}
	\begin{proof}
		As $A$ is an isometric isomorphism from $\calX_h$ to $\calY_1$ it is clear that $\left[ \begin{smallmatrix}
			A & 0 \\ 0 & I
		\end{smallmatrix}\right]z \in  C^1(\R_{\geq 0} ; \stack{\calY_1}{\calX_M}) \cap C(\R_{\geq 0} ; \dom \calB_M)$. Moreover, inserting $\left[ \begin{smallmatrix}
			A & 0 \\ 0 & I
		\end{smallmatrix}\right]z$ into \eqref{eq:pH_dynamics_2} yields the desired result by uniqueness of smooth solutions of \eqref{eq:ODEnode}.
	\end{proof}

	In other words, both formulations describe the same physical trajectory with the same boundary input and output signals. Observe that only the state coordinates differ (position-momentum $(x,M\dot x)$ versus strain-momentum $(Ax,M\dot x)$). Thus, the transformed node is not a different non-physical control system, but a different coordinate representation of the same one. It is the representation used, e.g., in \cite{kurula2015linear,skrepek2021linear}.

	\section{Applications}\label{sec:applications}
	In this part, we provide an application of the abstract theory developed in the previous section. In particular, we consider a~wave equation in Subsection~\ref{subsec:rod} and a Maxwell equation in Subsection~\ref{subsec:maxwell}.
	\subsection{$n$-dimensional wave equation with damping}\label{subsec:rod}
	As a first example, we consider a wave equation on a bounded Lipschitz domain $\Omega \subset \R^n$ modeled by the second-order system
	\begin{align}
		\rho(\zeta) \tfrac{\mathrm{d}^2}{\mathrm{d}t^2} x(t, \zeta) -\operatorname{div}( T(\zeta) \grad x(t, \zeta)) + a(\zeta)x(t, \zeta) + b(\zeta)\tfrac{\mathrm{d}}{\mathrm{d}t} x(t, \zeta) &= 0,  &t \in \R_+, \, &\zeta \in \Omega, \nonumber\\
		x(0, \zeta) &=x_0(\zeta) , && \zeta \in \Omega ,\nonumber\\
		\tfrac{\mathrm{d}}{\mathrm{d}t} x(0, \zeta) &=x_1(\zeta) , && \zeta \in \Omega    \nonumber
	\end{align}
	with displacement variable $x(t,\zeta)$. This system was recently treated in \cite{aigner2025well}, also in the context of boundary triplets.
	The above system is complemented by the boundary control and observation
	\begin{align}
	u(t)&=\tfrac12\big((1-P)\calR_{H^{1/2}(\partial\Omega)}\gamma_0\dot x(t)+(1+P)\,\gamma_\perp T\nabla x(t)\big),\label{eq:wave_bc} \\
	y(t)&=\tfrac12\big((1+P)\calR_{H^{1/2}(\partial\Omega)}\gamma_0\dot x(t)+(1-P)\,\gamma_\perp T\nabla x(t)\big), \notag
	\end{align}
	on $\partial\Omega$. Here, $\gamma_0\dot x$ denotes the boundary velocity and $\gamma_\perp T\nabla x=\nu\cdot T\nabla x|_{\partial\Omega}$ is the normal derivative (see Appendix~\ref{sec:app}). Moreover, $P\in L(H^{-1/2}(\partial\Omega))$ is a contraction specifying the type of boundary condition.

	We assume that the $a,b, \rho\in L^\infty(\Omega)$ are uniformly positive such that in particular, $\rho^{-1}\in L^\infty(\Omega)$. The same applies for the \emph{Young's modulus} $T\in L^\infty(\Omega)^{n\times n}$ for which we assume pointwise symmetry and uniform positive definiteness, such that $T^{-1}\in L^\infty(\Omega)$.
	Note that the space $H_0^{\operatorname{div}}(\Omega)$ is defined in~\eqref{eq:trace_kernels} in the Appendix \ref{sec:app}. We introduce the space $\calX \coloneqq  L^2(\Omega)$ and define the bounded multiplication operators
	\begin{align*}
		M: \calX \to \calX, \quad M x = \rho x, \qquad D: \calX \to \calX, \quad D x = b x,
	\end{align*}
	and the unbounded operator $S: \calX \supset \dom S \to \calX$
	\begin{align*}
		\dom S =\left\lbrace x \in H^1(\Omega) \, \middle| \, T \grad x \in H_0^{\operatorname{div}}(\Omega)   \right\rbrace , \qquad S x = -\operatorname{div}(T \grad x)  + a x.
	\end{align*}
	Here, $S$ is the unique self-adjoint operator that represents the form $h: \dom h \times \dom h \to \C$ in $\calX$.
    Moreover, as $a$ is uniformly positive, $\dom h =  H^1(\Omega)$ may also be endowed with the energy inner product
	\begin{align}\label{eq:innerprod}
		h(x,y) = \langle T \grad x,    \grad y \rangle_{L^2(\Omega ; \C^n)} + \langle ax , y \rangle_{L^2(\Omega)}
	\end{align}
	which we denote in the following by $\calX_h$.
	Hereby, we note that in~\cite{aigner2025well}, the authors choose an inner product on $\calX_h$ which includes a Dirichlet trace operator (instead of the second term involving $a$ in \eqref{eq:innerprod}).

	Using the above defined operators, the $n$-dimensional wave equation can be rewritten as the abstract second order system
	\begin{align*}
		M\Ddot{x}(t)+ D \dot{x}(t) + Sx(t)=0.
	\end{align*}
	At the same time, one should observe that the requirement of self-adjointness for $S$ restricts, from the outset, the admissible boundary conditions. Essentially, only functions with vanishing Neumann trace are permitted. To overcome this issue, we proceed as suggested in Section \ref{sec:main}. We emphasize that $S$ enters the construction only through its form domain, the energy space $\calX_h=\dom S^{1/2}=H^1(\Omega)$, which does not contain any boundary condition. The Neumann condition contained in $\dom S$ is not imposed on the solutions of the controlled system. Precisely, it is removed by passing to the maximal operator $\calA_{\max}$ below.
	We choose as state the variables position and momentum, i.e.
	\begin{equation*}
		z \coloneqq  \stack{
			x}{ M\dot{x} } .
	\end{equation*}
	Consequently, in view of \eqref{eq:energies}, we define the abstract potential energy $ H_p(z_1)= \tfrac{1}{2} h(z_1)$, which in this example coincides with the physical potential energy. In the spirit of \eqref{eq:h} we may also obtain the energy inner product by
	\begin{align*}
		h(x,y)=\langle Ax, Ay \rangle_{\calY}
	\end{align*}
	where
	$A: \calX = L^2(\Omega) \supset \calX_h = H^1(\Omega) \to \calY \eqqcolon \stack{ L^2(\Omega)^n}{  L^2(\Omega)}$ is given by
	\begin{equation*}
		A x = \begin{pmatrix}
			T^{1/2}\grad x \\ a^{1/2} x
		\end{pmatrix}.
	\end{equation*}
	We now investigate the precise representation of the adjoint $A^*$.
	\begin{lem}
		The adjoint of $A$ is given by $A^*: \calY \supset \dom A^* \to \calX$ where
		\begin{equation}
			\dom A^*= \stack{T^{-1/2}H_0^{\operatorname{div}}(\Omega)}{a^{-1/2}L^2(\Omega)}, \qquad A^*y = \begin{pmatrix}
				- \divergence T^{1/2} \  &  a^{1/2}
			\end{pmatrix}y.
		\end{equation}
	\end{lem}
	\begin{proof}
		Since $A: L^2(\Omega) \supset \dom A \to \calY$ can be seen as a column operator where only one entry is unbounded, we have $A^* = \nextto{(T^{1/2}\grad)^*}{a^{1/2}}$
		by \cite[Prop. 4.5]{moller2008adjoints}. It is well-known that the adjoint of $\grad : L^2(\Omega) \supset H^1(\Omega) \to L^2(\Omega)^n$ is $-\operatorname{div}: L^2(\Omega)^n \supset H_0^{\operatorname{div}}(\Omega) \to L^2(\Omega)$ which yields the claimed result.
	\end{proof}
	As a consequence,
	\begin{align*}
		A^*Ax = -\divergence (T \grad x) + a x , \qquad \dom A^*A = \left\lbrace x \in H^1(\Omega) \, \middle| \, \stack{ T \grad x}{a x} \in \stack{H_0^{\operatorname{div}}(\Omega)}{L^2(\Omega)}  \right\rbrace
	\end{align*}
	is the self-adjoint representing operator of the form $h$ corresponding to the potential energy. The operator $A$ restricts to a bounded isometry $A:\calX_h \to \calY_1$. Then, the operator $A^*A$ uniquely extends to the Riesz isomorphism $\calR_{\calX_h} : \calX_h \to \calX_h^*$.

	To introduce the state space corresponding to the momentum variable, we let $\calX_M=L^2(\Omega)$ endowed with the equivalent inner product
	\begin{align*}
		\langle x, y \rangle_{\calX_M}=\langle M^{-1} x, y \rangle_\calX.
	\end{align*}
	We have the kinetic energy \begin{align*}
		H_k(z_2)= \frac{1}{2} \|z_2\|_{\calX_M}^2 = \frac{1}{2} \|\rho^{-1/2}z_2\|_{\calX}^2.
	\end{align*}
	Moreover, to define the dual pair satisfying \eqref{eq:dual_pair_ass}, we see that the operator $B: \calY \supset \dom B \to \calX$ where
	\begin{align*}
		\dom B = \stack{T^{-1/2}H^{\operatorname{div}}(\Omega)}{a^{-1/2}L^2(\Omega)}, \qquad By=\begin{pmatrix}
			\divergence T^{1/2} \ &  -a^{1/2}
		\end{pmatrix}y
	\end{align*}
	extends $-A^*$. Hence $(A^*,-B^*)$ forms a dual pair. It follows that $\calA_\mathrm{max}: \stack{\calX_h}{\calX_M} \supset \dom{\calA_\mathrm{max}}\to \stack{\calX_h}{\calX_M}$ where
	\begin{align*}
		\calA_\mathrm{max}\roundstack{z_1}{z_2}= \begin{pmatrix}
			\rho^{-1}z_2 \\ \divergence T \grad z_1  - az_1
		\end{pmatrix}, \qquad  \dom \calA_\mathrm{max}= \left\lbrace \roundstack{z_1}{z_2} \in \stack{\calX_h}{M \calX_h}  \, \middle| \,  T\grad z_1 \in H^{\operatorname{div}}(\Omega) \right\rbrace
	\end{align*}
	is closed, densely defined and the adjoint of a skew-symmetric operator by Lemma \ref{lem:A_dense}. Furthermore, the integration by parts formula (see \eqref{eq:int_by_parts} in Appendix~\ref{sec:app}) immediately yields that there is a boundary triplet for the dual pair $(A^*, -B^*)$ as proven in the following.
	\begin{lem}
		The triplet $(\calG_1, \Lambda, \Pi)= (H^{1/2}(\partial \Omega), \gamma_0, \nextto{- \gamma_\perp  T^{1/2}}{0})$ defines a boundary triplet for the dual pair $(A^*,-B^*)$.
	\end{lem}
	\begin{proof}
		We have to check that the linear operators
		\begin{align*}
			\Lambda=\Lambda_1= \gamma_0 : H^1(\Omega) \to H^{1/2}(\partial \Omega), \quad \Pi=\Pi_1 = \nextto{- \gamma_\perp  T^{1/2} }{0}: \stack{T^{-1/2}H^{\operatorname{div}}(\Omega)}{ a^{-1/2}L^2(\Omega)} \to H^{-1/2}(\partial \Omega)
		\end{align*}
		satisfy
		\begin{enumerate}
			\item[i)] $\ran \Lambda = H^{1/2}(\partial \Omega)$, $ \ran \Pi = H^{-1/2}(\partial \Omega)$;
			\item[ii)] the abstract Green identity \begin{align}
				-\langle By,x \rangle_\calX - \langle y, Ax \rangle_\calY = \langle \Pi y, \Lambda x \rangle_{\calG_1^*, \calG_1}
			\end{align}
			holds for all $y \in \dom B, x \in \dom A$.
		\end{enumerate}
		The abstract green identity follows from
		\begin{align*}
			- \langle \divergence  T^{1/2}  y_1 - a^{1/2}y_2 , x\rangle_{\calX} - \langle \roundstack{y_1}{y_2} , \roundstack{T^{1/2} \grad x}{a^{1/2}x} \rangle_{\calY}
			&= - (\langle \divergence  T^{1/2} y_1 , x\rangle_{L^2(\Omega)} + \langle y_1 ,  T^{1/2} \grad x \rangle_{L^2(\Omega)^n} ) \\
			&= \langle \nextto{-\gamma_\perp  T^{1/2}}{0} \roundstack{y_1}{y_2} , \gamma_0 x\rangle_{H^{-1/2}(\partial \Omega),H^{1/2}(\partial \Omega)}
		\end{align*}
		for all $\roundstack{y_1}{y_2} \in \stack{ T^{-1/2} H^{\operatorname{div}}(\Omega)}{ a^{-1/2}L^2(\Omega)}, x\in H^1(\Omega)$, see \eqref{eq:int_by_parts} of Appendix~\ref{sec:app}. The surjectivity is immediate from Theorem \ref{thm:traces} in Appendix~\ref{sec:app}.
	\end{proof}

	We are now in a position to apply Theorem~\ref{thm:bdd_triplet_for_sec_order}, which in turn allows us to provide an explicit characterization of $\dom \calA_\mathrm{max}^*$. For this purpose, we introduce the space
	\begin{align}
		H_\mathrm{N}^1(\Omega) \coloneqq  \left\lbrace z_1 \in H^1(\Omega) \, \middle| \, T \grad z_1 \in H^{\operatorname{div}} (\Omega), \, \gamma_\perp T \grad z_1 = 0  \right\rbrace
	\end{align}
	and we note that the equality $\gamma_\perp T \grad z_1 = 0$ is meant in $H^{-1/2}(\partial \Omega)$.
	\begin{lem}
		Define $\Gamma \coloneqq  \stack{\Gamma_0}{\Gamma_1}: \dom \calA_\mathrm{max}\to H^{1/2}(\partial \Omega) \times H^{-1/2}(\partial \Omega)$ by
		\begin{align*}
			\Gamma z = \begin{bmatrix}
				\Gamma_0 \\ \Gamma_1
			\end{bmatrix} z = \begin{bmatrix}
				0 &  \gamma_0 \\
				\gamma_\perp T \grad &  0
			\end{bmatrix} z
		\end{align*}
		Then, $(H^{1/2}(\partial \Omega), \Gamma_0, \Gamma_1)$ is a boundary triplet for $\calA_0^*$. In addition,
		\begin{align}\label{eq:min_dom1}
			\dom \calA_0=\dom \calA_\mathrm{max}^*
			= H_\mathrm{N}^1(\Omega) \times H_0^1(\Omega).
		\end{align}
	\end{lem}
	\begin{proof}
		The first assertion is true by Theorem \ref{thm:bdd_triplet_for_sec_order}. The second assertion \eqref{eq:min_dom1} follows as the domain of the minimal operator $\calA_0$ can be recovered by restricting $\dom \calA_\mathrm{max}$ to $\ker \Gamma$ by Lemma \ref{lem:cont_bdd_map}.
	\end{proof}
	We now introduce the  boundary node associated with the $n$-dimensional wave equation. To this end, let $P \in L(H^{-1/2}(\partial \Omega))$ be a contraction and define
	\begin{align*}
		N = \left( \roundstack{z_1}{z_2} \mapsto \begin{bsmallmatrix}
			\tfrac{1}{2}((1-P)\calR_{H^{1/2}(\partial \Omega)}\gamma_0 \rho^{-1} z_2 + (1+P) \gamma_\perp  T \grad z_1 ) \\ \begin{pmatrix}
				\rho^{-1}z_2 \\ \divergence (T \grad z_1 ) - az_1 - b\rho^{-1}z_2
			\end{pmatrix} \\ \tfrac{1}{2}((1+P)\calR_{H^{1/2}(\partial \Omega)}\gamma_0 \rho^{-1} z_2 + (1-P) \gamma_\perp T  \grad z_1 )
		\end{bsmallmatrix}, \begin{bsmallmatrix}
			H^{-1/2}(\partial \Omega) \vphantom{\tfrac{1}{2}((1-P)\calR_{H^{1/2}(\partial \Omega)}\gamma_0 \rho^{-1} z_2 + (1+P) \gamma_\perp  T \grad z_1 )} \\ \stack{H^1(\Omega)}{L^2(\Omega)} \vphantom{\begin{pmatrix}
				\rho^{-1}z_2 \\ \divergence (T \grad z_1 ) - az_1 - b\rho^{-1}z_2
			\end{pmatrix}} \\ H^{-1/2}(\partial \Omega) \vphantom{\tfrac{1}{2}((1-P)\calR_{H^{1/2}(\partial \Omega)}\gamma_0 \rho^{-1} z_2 + (1+P) \gamma_\perp  T \grad z_1 )}
		\end{bsmallmatrix}\right)
	\end{align*}
	with
	\begin{equation*}
		\dom N = \left\lbrace \roundstack{z_1}{z_2} \in \stack{H^1(\Omega)}{MH^1(\Omega)} \, \middle| \,
		T \grad z_1 \in H^{\operatorname{div}}(\Omega)
		\right\rbrace .
	\end{equation*}
	Then, $N$ defines an impedance passive internally well-posed boundary node by Theorem \ref{thm:main2}. Moreover, the main operator of $N$ is the generator of a unitary semigroup on $\stack{H^1(\Omega)}{L^2(\Omega)}$ if we choose a unitary operator $P$ and $b \equiv 0$. In particular, the choice $P=1$ leads to the situation as described above:
	\begin{align*}
		\dom S = \ker \gamma_\perp T \grad = H_\mathrm{N}^1(\Omega) .
	\end{align*}
    \begin{rem}\label{rem:damping_bc}

	As discussed in the introduction, the choice $P=1$ correpsonds to Neumann boundary control and velocity observation, whereas $P=-1$ is control of the boundary velocity.

    However, let $k\in L^\infty(\partial\Omega)$ with $k\ge0$ almost everywhere, and let $K\colon H^{1/2}(\partial\Omega)\to H^{-1/2}(\partial\Omega)$, $Kv=kv \in L^2(\partial\Omega)\subset H^{-1/2}(\partial\Omega)$. Then $K$ is bounded so that $\calR+K$ is boundedly invertible. Define the operator $P\coloneqq (\calR-K)(\calR+K)^{-1}\in L(H^{-1/2}(\partial\Omega))$ and let $g\coloneqq  (\calR+K)v$ for some $v \in H^{1/2}(\partial\Omega)$. It holds that
    \begin{equation*}
    \begin{aligned}
    \|Pg\|_{H^{-1/2}(\partial\Omega)}^2-\|g\|_{H^{-1/2}(\partial\Omega)}^2
    &=\Re\big\langle Pg-g,\;Pg+g\big\rangle_{H^{-1/2}(\partial\Omega)}
    =-\Re\big\langle 2Kv,\;2\calR v\big\rangle_{H^{-1/2}(\partial\Omega)}\\
    &=-4\,\Re\langle Kv,v\rangle_{H^{-1/2}(\partial\Omega),H^{1/2}(\partial\Omega)}
    =-4\int_{\partial\Omega}k|v|^2\,\mathrm d\sigma\le 0 .
    \end{aligned}
    \end{equation*}
    Using $1+P=2\calR(\calR+K)^{-1}$, $1-P=2K(\calR+K)^{-1}$ and $K(\calR+K)^{-1}\calR=\calR(\calR+K)^{-1}K$, the input equation \eqref{eq:wave_bc} becomes
	\[
	u=\calR(\calR+K)^{-1}\big(\gamma_\perp T\nabla x+k\,\gamma_0\dot x\big).
	\]
	Since $\calR(\calR+K)^{-1}$ is an isomorphism on $H^{-1/2}(\partial\Omega)$, this is, up to a bounded and boundedly invertible transformation of the input, the boundary control
	\[
	\nu\cdot T\nabla x+k\,\dot x=\tilde u\qquad\text{on }\partial\Omega.
	\]
	The resulting boundary condition in the zero input case is $\nu\cdot T\nabla x+k\,\dot x=0$. Note that this is the classical dissipative boundary condition, which is obtained from Neumann control by the feedback $u=-k\dot x$ and is used for uniform boundary stabilization of the wave equation, see, e.g., \cite[Eq. (1.2), (1.6)]{komornik1990direct} and \cite{lagnese1983decay}. Theorem~\ref{thm:main2} shows that the corresponding open system is impedance passive and internally well-posed. Observe that $P$ is isometric if and only if $k=0$.
	\end{rem}
	In addition, recall that
	\begin{align*}
		\calY_1 = \ran A= \left\lbrace \roundstack{T^{1/2}\grad z_1}{a^{1/2}z_1} \, \middle| \, z_1 \in H^1(\Omega) \right\rbrace
	\end{align*}
	is a Hilbert space endowed with he usual inner product on $\calY$. Furthermore,
	\begin{align*}
		\dom B_{\calY_1}=  \left\lbrace \roundstack{w_1}{w_2} \in \stack{T^{-1/2}H^{\operatorname{div}}(\Omega)}{ a^{-1/2}L^2(\Omega)} \, \middle| \, \exists \, z_1 \in H^1(\Omega) : \roundstack{w_1}{w_2} =\roundstack{T^{1/2}\grad z_1}{a^{1/2}z_1}  \right\rbrace  .
	\end{align*}
	According to Theorem \ref{thm:trafo} we define
	\begin{align*}
		L = \left( \roundstack{\stack{w_1}{w_2}}{w_3} \mapsto \begin{bsmallmatrix}
			\tfrac{1}{2}((1-P)\calR_{H^{1/2}(\partial \Omega)}\gamma_0 \rho^{-1} w_3 + (1+P) \gamma_\perp  T^{1/2} w_1 ) \\ \begin{pmatrix}
				T^{1/2} \grad \rho^{-1}w_3 \\ a^{1/2} \rho^{-1}w_3  \\ \divergence T^{1/2} w_1 - a^{1/2} \rho^{-1}w_2 - b\rho^{-1}w_2
			\end{pmatrix} \\\tfrac{1}{2}((1+P)\calR_{H^{1/2}(\partial \Omega)}\gamma_0 \rho^{-1} w_3 + (1-P) \gamma_\perp  T^{1/2} w_1 )
		\end{bsmallmatrix}, \begin{bsmallmatrix}
			H^{-1/2}(\partial \Omega) \vphantom{\tfrac{1}{2}((1-P)\calR_{H^{1/2}(\partial \Omega)}\gamma_0 \rho^{-1} w_3 + (1+P) \gamma_\perp  T^{1/2} w_1 )} \\ \stack{\calY_1}{L^2(\Omega)} \vphantom{\begin{pmatrix}
				T^{1/2} \grad \rho^{-1}w_3 \\ a^{1/2} \rho^{-1}w_3  \\ \divergence T^{1/2} w_1 - a^{1/2} \rho^{-1}w_2 - b\rho^{-1}w_2
			\end{pmatrix}} \\ H^{-1/2}(\partial \Omega) \vphantom{\tfrac{1}{2}((1-P)\calR_{H^{1/2}(\partial \Omega)}\gamma_0 \rho^{-1} w_3 + (1+P) \gamma_\perp  T^{1/2} w_1 )}
		\end{bsmallmatrix}\right) , \qquad  \dom L = \stack{\dom B_{\calY_1}}{H^1(\Omega)} .
	\end{align*}
	Then, also $L$ defines an impedance passive internally well-posed boundary node.

	\subsection{Maxwell equation}\label{subsec:maxwell}
	Let $\Omega\subset\R^3$ be a bounded and simply connected Lipschitz domain with connected boundary. We consider the \textit{Maxwell equation} in the vacuum, i.e.,
	\begin{align}
		\tfrac{\mathrm{d}^2}{\mathrm{d}t^2} \mathbf{B}(t, \zeta) + c^2 \curl \curl \mathbf{B}(t, \zeta)    &=0, &t \in \R_+, \, &\zeta \in \Omega, \nonumber\\
\operatorname{div}\mathbf{B}(t, \zeta)   &=0, &t \in \R_+, \, &\zeta \in \Omega, \nonumber\\
		\mathbf{B}(0, \zeta) &=\mathbf{B}_0(\zeta) , && \zeta \in \Omega , \label{eq:Max}\\
		\tfrac{\mathrm{d}}{\mathrm{d}t} \mathbf{B}(0, \zeta) &=\mathbf{B}_1(\zeta) , && \zeta \in \Omega. \nonumber
	\end{align}
Hereby
$\mathbf{B}$ stands for the magnetic flux density, and $c>0$ is the speed of light. Divergence-freeness of $\mathbf{B}$ will be encoded in the space, which is chosen as $\calX \coloneqq  L^2(\operatorname{div}\!=\!0;\Omega)^3$, see Appendix~\ref{sec:app}. Since no operator occurs in front of the second temporal derivative, we set
$M=I_\calX$.
The form associated to the potential energy $h: \dom h \times \dom h \to \C$ with $\dom h = H(\curl,\operatorname{div}\!=\!0; \Omega)$ via
	\begin{align*}
		h(x,y) = c^2\,\langle  \curl x,    \curl y \rangle_{L^2(\Omega)^3}.
	\end{align*}
	Closedness and coercivity of $h$ follows from Lemma~\ref{lem:curlcoerc} together with the fact that $H(\curl,\operatorname{div}\!=\!0;\Omega)$ is a~Hilbert space. In particular, we have that
        $\calX_h =(H(\curl,\operatorname{div}\!=\!0;\Omega), h)$ is a Hilbert space. Moreover, $h$ is represented by the operator
	\begin{align*}
		A: L^2(\operatorname{div}\!=\!0;\Omega) \supset H(\curl,\operatorname{div}\!=\!0; \Omega) \to L^2(\operatorname{div}\!=\!0,\Omega) \qquad  A x = c\,\curl x
	\end{align*}
	in the sense that $S=A^*A$ satisfies $ h(x,y) = \langle Sx,y  \rangle_{L^2(\Omega)^3}$ for all $x \in \dom A^*A = \dom S$ and $y \in \dom h$. The adjoint is given by
\begin{align}\label{eq:Maxadjoint}
		A^*: {L^2(\operatorname{div}\!=\!0;\Omega)} \supset {H_0(\curl,\operatorname{div}\!=\!0; \Omega)} \to L^2(\operatorname{div}\!=\!0;\Omega), \qquad  A^* y = c\,\curl y.
	\end{align}
	We abbreviate $\calY \coloneqq  {L^2(\operatorname{div}\!=\!0;\Omega)}$. By defining $x(t)=\mathbf{B}(t, \cdot)$, \eqref{eq:Max} may be rewritten as the abstract second-order system
	\begin{align*}
		\Ddot{x}(t)+A^*A x(t)=0,
	\end{align*}
	which, through the domain of the adjoint operator~\eqref{eq:Maxadjoint} implicitly imposes the boundary condition $\gamma_\times c \curl \mathbf{B}=\gamma_\pi c \curl \mathbf{B}=0$. To overcome this issue and to be able to impose a wider range of boundary conditions, we introduce the operator
	\begin{align*}
		B:  \calY \supset {H(\curl,\operatorname{div}\!=\!0; \Omega)} \to L^2(\operatorname{div}\!=\!0;\Omega), \qquad  B y =-c\curl y.
	\end{align*}
    In the following, we will denote by $\gamma_\times$ and $\gamma_\pi$ the tangential boundary trace and tangential trace projection as introduced in Theorem~\ref{thm:tangtraces} and Lemma~\ref{lem:div0trace} in Appendix~\ref{sec:app}. Moreover, $V_\pi(\partial \Omega)$ and $V_\times(\partial \Omega)$ denote the range spaces onto which these traces map surjectively. In view of Theorem~\ref{thm:tangtraces} and \eqref{eq:int_by_parts_curl} in Appendix~\ref{sec:app}, it is immediate to verify that that $(A^*, -B^*)$ forms a dual pair, i.e. $A^* \subset -B$.
	\begin{lem}
		The triple $(\calG_1, \Lambda, \Pi)= (V_\pi(\partial \Omega), \gamma_\pi,c {\gamma_\times})$ defines a boundary triplet for the dual pair $(A^*,-B^*)$.
	\end{lem}
	\begin{proof}
	Surjectivity of $\gamma_\pi$ and $\gamma_\times$ has been proven in Lemma~\ref{lem:div0trace}, whereas the abstract Green identity follows from \eqref{eq:int_by_parts_curl}.
	\end{proof}
	Now we can introduce the operator corresponding to the first order formulation, i.e., $\calA_\mathrm{max}: \stack{\calX_h}{\calX} \supset  \dom \calA_\mathrm{max} \to \stack{\calX_h}{\calX}$ via
	\begin{align*}
		 \dom \calA_\mathrm{max}= \left\lbrace \roundstack{z_1}{z_2} \in \stack{H(\operatorname{curl},\operatorname{div}\!=\!0;\Omega) }{H(\operatorname{curl},\operatorname{div}\!=\!0;\Omega) }  \, \middle| \, c \curl z_1 \in H(\operatorname{curl},\operatorname{div}\!=\!0;\Omega)  \right\rbrace .
	\end{align*}
    and action
    \begin{equation*}
        \calA_\mathrm{max}\roundstack{z_1}{z_2} = \begin{pmatrix}
			z_2 \\ - c^2\,\curl \curl z_1
		\end{pmatrix}.
    \end{equation*}
	\begin{lem}
		The triplet $(V_\pi(\partial \Omega),  \nextto{0}{\gamma_\pi},\nextto{c {\gamma_\times \curl}}{0})$ is a boundary triplet for $\calA_\mathrm{max}$.
	\end{lem}
	\begin{proof}
		The claim follows from Theorem \ref{thm:bdd_triplet_for_sec_order}.
	\end{proof}
	We now introduce the associated boundary node using Theorem \ref{thm:main2}. Let $P \in L(V_\times(\partial \Omega))$ be a contraction and define the colligation
	\begin{align*}
		N = \left( \roundstack{z_1}{z_2} \mapsto \begin{bsmallmatrix}
			\tfrac{1}{2}((1-P)\calR_{V_\times(\partial \Omega)}^{-1} \gamma_\pi z_2 + (1+P) c \gamma_\times  \curl z_1 ) \\ \begin{pmatrix}
				& I \\ - c^2 \curl \curl  &
			\end{pmatrix}\begin{pmatrix}
				z_1 \\ z_2
			\end{pmatrix}
			\\ \tfrac{1}{2}((1+P)\calR_{V_\times(\partial \Omega)}^{-1} \gamma_\pi z_2 + (1-P) c \gamma_\times   \curl z_1 )
		\end{bsmallmatrix}, \begin{bsmallmatrix}
			V_\times(\partial \Omega) \vphantom{\tfrac{1}{2}((1-P)\calR_{V_\times(\partial \Omega)}^{-1} \gamma_\pi z_2 + (1+P) c \gamma_\times  \curl z_1 )} \\ \stack{H(\operatorname{curl},\operatorname{div}\!=\!0;\Omega)}{L^2(\Omega)^3} \vphantom{\begin{pmatrix}
				& I \\ - c^2 \curl \curl  &
			\end{pmatrix}\begin{pmatrix}
				z_1 \\ z_2
			\end{pmatrix} } \\ V_\times(\partial \Omega) \vphantom{\tfrac{1}{2}((1-P)\calR_{V_\times(\partial \Omega)}^{-1} \gamma_\pi z_2 + (1+P) c \gamma_\times  \curl z_1 )}
		\end{bsmallmatrix}\right)
	\end{align*}
	with
	\begin{equation*}
		\dom N = \dom \calA_\mathrm{max}.
	\end{equation*}
	Hence following Theorem \ref{thm:main2}, $N$ defines an impedance passive internally well-posed boundary node. Recall that $\calY_1 = c \curl H(\operatorname{curl},\operatorname{div}\!=\!0;\Omega)$ and as the operator $A$ uniquely restricts to a bounded isometry $A: H(\operatorname{curl},\operatorname{div}\!=\!0;\Omega) \to \calY_1$ we introduce the similarity transform
	\begin{align*}
		L = \left( \roundstack{w_1}{w_2} \mapsto \begin{bsmallmatrix}
			\tfrac{1}{2}((1-P)\calR_{V_\times(\partial \Omega)}^{-1} \gamma_\pi w_2 + (1+P) c \gamma_\times w_1 ) \\ \begin{pmatrix}
				& c\curl \\ -c\curl & \end{pmatrix} \begin{pmatrix}
				w_1 \\ w_2
			\end{pmatrix} \\
			\tfrac{1}{2}((1+P)\calR_{V_\times(\partial \Omega)}^{-1} \gamma_\pi w_2 + (1-P) c \gamma_\times w_1 )
		\end{bsmallmatrix}, \begin{bsmallmatrix}
			V_\times(\partial \Omega) \vphantom{\tfrac{1}{2}((1-P)\calR_{V_\times(\partial \Omega)}^{-1} \gamma_\pi w_2 + (1+P) c \gamma_\times w_1 )} \\ \stack{\calY_1}{L^2(\Omega)^3} \vphantom{\begin{pmatrix}
				& c\curl \\ -c\curl & \end{pmatrix} \begin{pmatrix}
				w_1 \\ w_2
			\end{pmatrix}} \\ V_\times(\partial \Omega) \vphantom{\tfrac{1}{2}((1-P)\calR_{V_\times(\partial \Omega)}^{-1} \gamma_\pi w_2 + (1+P) c \gamma_\times w_1 )}
		\end{bsmallmatrix}\right), \qquad  \dom L = \stack{\dom B_{\calY_1}}{H(\operatorname{curl},\operatorname{div}\!=\!0;\Omega)},
	\end{align*}
	where
	\begin{align*}
		\dom B_{\calY_1}=  \left\lbrace w_1 \in H(\operatorname{curl},\operatorname{div}\!=\!0;\Omega) \, \middle| \, \exists \, z_1 \in H(\operatorname{curl},\operatorname{div}\!=\!0;\Omega) : w_1=c\curl z_1 \right\rbrace  .
	\end{align*}
	Consequently, $L$ defines an impedance passive internally well-posed boundary node. Note that the input boundary map of the colligation $L$ is very similar to the \textit{Leontovich boundary condition} in \cite[Thm. 2.6]{eller2022m}.

	\section{Conclusion}\label{sec:conclusion}
	In this work, we consider second-order evolution equations by means of their first-order in time representation. We constructed a boundary triplet for the associated block operator matrix and, based on the abstract boundary mappings, defined natural input and output operators. This construction provided a systematic way to open the previously closed system to boundary control interaction. Furthermore, we introduced a similarity transformation for the resulting boundary control system, derived from a factorization of the stiffness operator motivated by standard PDE models. The applicability of the proposed framework was demonstrated through two examples.

\section*{Declarations}

 \paragraph*{Data Availability Statement: } \;
There are no data equipped with this article. Data sharing is therefore not applicable to this article.

 \paragraph*{Underlying and related material:}\;
 No underlying or related material.

 \paragraph*{Author contributions:}\;
 The authors Till Preuster, Timo Reis and Manuel Schaller contributed equally.

\paragraph*{Competing interests:}\;
No competing interests to declare.

	\appendix

	\section{Sobolev spaces and trace theorems}\label{sec:app}
	Let $\Omega \subset \R^n$ be a {bounded Lipschitz domain}.
	As usual, we denote the weak divergence, gradient, Laplacian, and curl by $\operatorname{div}$, $\nabla$, $\Delta$, and $\operatorname{curl}$, respectively.
	Clearly, $H^1(\Omega)$ denotes the space of $L^2$-functions whose weak gradient belongs to $L^2$.
	In analogy, we consider
	\begin{align*}
		H^{\operatorname{div}} (\Omega) &= \left\lbrace x \in L^2(\Omega)^n \, \middle| \, \operatorname{div} x \in L^2(\Omega) \right\rbrace, &
		H(\operatorname{curl};\Omega) &= \left\lbrace x \in L^2(\Omega)^3 \, \middle| \, \operatorname{curl} x \in L^2(\Omega)^3 \right\rbrace .
	\end{align*}
	These spaces become Hilbert spaces when equipped with the corresponding graph norms.

	The space $H_0^1(\Omega)$ is usually defined as the closure of $C_0^\infty(\Omega)$ with respect to the $H^1(\Omega)$-norm.
	In analogy, we introduce
	\begin{align}\label{eq:trace_kernels}
		H_0^{\operatorname{div}}(\Omega)\coloneqq  \overline{C_0^\infty(\Omega)^n}^{H^{\operatorname{div}}(\Omega)}, \qquad  H_0^1(\Omega)\coloneqq  \overline{C_0^\infty(\Omega)}^{H^1(\Omega)}, \qquad  H_0(\operatorname{curl};\Omega)\coloneqq  \overline{C_0^\infty(\Omega)}^{H(\operatorname{curl};\Omega)}.
	\end{align}
	We refer the reader to \cite{adams2003sobolev} for the definition of fractional Sobolev spaces. Recall that, the space $H^{-1/2}(\partial \Omega)$ is the dual space of $H^{1/2}(\partial \Omega)$ with respect to the pivot space $L^2(\partial \Omega)$.
	Let $\nu \in L^\infty(\partial\Omega)^3$ denote the outward unit normal vector, which is well-defined almost everywhere with respect to the surface measure on~$\partial\Omega$, since $\Omega$ is a Lipschitz domain. We recall an important results from trace theory.
	\begin{thm}\label{thm:traces}
		The following is true for the bounded Lipschitz domain $\Omega \subset \R^n$:
		\begin{enumerate}[(a)]
			\item The boundary trace mapping $x \mapsto \left. x \right|_{\partial \Omega}: C^1(\overline{\Omega})\to  L^2(\partial\Omega)$ has a unique bounded and surjective extension $\gamma_0:H^1(\Omega) \to H^{1/2}(\partial \Omega)$. Moreover, $H_0^1(\Omega)=\ker \gamma_0$.
			\item The normal trace mapping $y \mapsto \nu^\top \gamma_0 y: H^1(\Omega)^n \to L^2(\partial \Omega)$ has a unique bounded and surjective extension $\gamma_{\perp}: H^{\operatorname{div}}(\Omega) \to H^{-1/2}(\partial \Omega)$.
		\end{enumerate}	\end{thm}
	The above result is scattered across the literature: the surjectivity of the boundary trace mapping is established in \cite[Thm.~1.5.1.2]{grisvard_elliptic_1985}, the identity $H_0^1(\Omega) = \ker \gamma_0$ is proved in \cite[Thm.~5.37]{adams2003sobolev}, and the statement concerning the normal trace mapping follows from \cite[Lem.~20.2]{Tartar2007Sobolev}.
	The latter theorem also shows that the following integration by parts formula holds:
	\begin{align}\label{eq:int_by_parts}
		\forall\,x \in H^{\operatorname{div}}(\Omega), y \in H^1(\Omega):\quad
		\langle \operatorname{div} x, y \rangle_{L^2(\Omega)} + \langle x, \grad y \rangle_{L^2(\Omega)^n}=\langle \gamma_\perp x, \gamma_0 y \rangle_{H^{-1/2}(\partial \Omega),H^{1/2}(\partial \Omega)}.
	\end{align}
	Trace theorems and an integration by parts formula for $H(\operatorname{curl}; \Omega)$ are somewhat more involved for two reasons. First, only the tangential component of the boundary trace must be considered, rather than the full trace. Second, it turns out that the resulting trace spaces are neither subsets nor supersets of $L^2(\partial \Omega, \R^3)$. The tangential trace space associated with $H(\operatorname{curl}; \Omega)$ is studied in \cite{BuffaCostabelSheen2002}.
	Although our approach is more abstract, it is somewhat simpler, as we do not explicitly construct these spaces.
	For a function $x \in C^1(\overline{\Omega})$, we define its \emph{tangential trace} $x_\times \in L^2(\partial\Omega)^3$ by
	\begin{equation}
		x_\times(\xi) = \nu(\xi) \times x(\xi) \quad \text{for almost all } \xi \in \partial\Omega.
		\label{eq:xtimes}
	\end{equation}
	Further, the {\em tangential trace projection}
	of $x\in C^1(\overline{\Omega})$ is $x_\pi\in L^2(\partial\Omega)^3$ with
	\begin{equation}
		x_\pi(\xi)=(\nu(\xi)\times x(\xi))\times \nu(\xi)\text{ for almost all }\xi\in \partial\Omega.\label{eq:xpi}\end{equation}
	Both sets of tangential traces and tangential trace projections of smooth functions are not dense in $L^2(\partial\Omega)^3$, because
	such functions are almost everywhere pointwise in the tangential space of $\partial\Omega$. This motivates to introduce
	\[
	L^2_\tau(\partial\Omega) \coloneqq
	\left\{
	w \in L^2(\partial\Omega)^3 \;\middle|\; \nu^\top w = 0 \text{ on } \partial\Omega
	\right\},
	\]
	which is a~closed subspace of $L^2(\partial\Omega)^3$, and therefore a~Hilbert space endowed with the standard inner product in $L^2(\partial\Omega)^3$.\\
	Next we introduce the trace theorem for $H(\operatorname{curl}; \Omega)$.
	\begin{thm}\label{thm:tangtraces}
		Let $\Omega \subset \R^n$ be a~bounded Lipschitz domain.
		There exist Hilbert spaces $V_\times(\partial\Omega)$, $V_\pi(\partial\Omega)$ with the following properties:
		\begin{enumerate}[(a)]
		\item\label{thm:tangtracesa} $L^2_\tau(\partial\Omega)\cap V_\times(\partial\Omega)\cap V_\pi(\partial\Omega)$ is dense in all  $L^2_\tau(\partial\Omega)$, $V_\times(\partial\Omega)$ and $V_\pi(\partial\Omega)$.
			\item\label{thm:tangtracesb} $V_\pi(\partial\Omega)$ is the dual space of $V_\times (\partial\Omega)$, and the duality product fulfills
			\[\langle z_1,z_2\rangle_{V_\times,V_\pi}=\langle z_1,z_2\rangle_{L^2(\partial\Omega)^3}\text{ for all }z_1\in L^2_\tau(\partial\Omega)\cap V_\times(\partial\Omega),\,z_2\in L^2_\tau(\partial\Omega)\cap V_\pi(\partial\Omega).\]
			\item\label{thm:tangtracesc} The tangential boundary trace mapping $x \mapsto  x_\times: C^1(\overline{\Omega}) \to L^2_\tau(\partial \Omega)$ (with $x_\times$ as in \eqref{eq:xtimes}) has a unique bounded and surjective extension $\gamma_\times:H(\operatorname{curl}; \Omega) \to V_\times(\partial \Omega)$.
			\item\label{thm:tangtracesd} The tangential trace projection mapping $x \mapsto  x_\pi: C^1(\overline{\Omega}) \to L^2_\tau(\partial \Omega)$ (with $x_\pi$ as in \eqref{eq:xpi}) has a unique bounded and surjective extension $\gamma_\pi:H(\operatorname{curl}; \Omega) \to V_\pi(\partial \Omega)$.
		\end{enumerate}
		Moreover,
		\begin{equation}
			H_0(\operatorname{curl};\Omega)=\ker \gamma_\times=\ker \gamma_\pi.\label{eq:kergamma}
		\end{equation}
	\end{thm}
	The construction of the spaces $V_\times(\partial\Omega)$ and $V_\pi(\partial\Omega)$, as well as the proofs of their properties stated in the above theorem, were carried out in \cite[Section 4.4]{PhilReis23}, motivated by the results of \cite{skrepek2021linear}. Furthermore, the relation \eqref{eq:kergamma} follows from \cite[Theorem I.2.6]{GiraultRaviart1986}.\\
	Although $V_\pi(\partial\Omega) = V_\times(\partial\Omega)^*$, the triple $(V_\times(\partial\Omega), L^2_\tau(\partial\Omega), V_\pi(\partial\Omega))$ does not form a Gelfand triple, as there is no inclusion relation between the spaces. Instead, it constitutes a quasi Gelfand triple in the sense of \cite{Skrepek2024}.
	We also note that there is an alternative approach to constructing trace spaces for $H(\operatorname{curl}; \Omega)$ based on the theory of {\em Hilbert complexes}; see \cite{HiptmairPaulySchulz2023}, where suitable trace spaces were derived in an abstract manner by considering the quotient space $H(\operatorname{curl}; \Omega) / H_0(\operatorname{curl}; \Omega)$.\\
	Next we show that, by using the above presented spaces, the curl operator fulfills an integration by parts formula, see e.g.\ \cite[Section 4.4]{PhilReis23}. Namely, for any bounded Lipschitz domain $\Omega\subset\R^3$, it holds that
	\begin{align}\label{eq:int_by_parts_curl}
		\forall\,x,y \in H({\operatorname{curl}};\Omega):\quad        \langle \operatorname{curl} x, y \rangle_{L^2(\Omega)^3} - \langle x, \operatorname{curl} y \rangle_{L^2(\Omega)^3}=\langle \gamma_\times x, \gamma_\pi y \rangle_{V_\times(\partial \Omega),V_\pi(\partial \Omega)}.
	\end{align}
Next we place our attention to spaces of divergence-free functions. To this end we consider the spaces
\begin{align*}
L^2(\operatorname{div}\!=\!0;\Omega) &\coloneqq \setdef{x\in H^{\operatorname{div}}(\Omega)}{\operatorname{div} x=0},\\
H(\operatorname{curl},\operatorname{div}\!=\!0;\Omega) &\coloneqq H(\operatorname{curl};\Omega)\phantom{_0}\cap L^2(\operatorname{div}\!=\!0;\Omega),\\
H_0(\operatorname{curl},\operatorname{div}\!=\!0;\Omega) &\coloneqq H_0(\operatorname{curl};\Omega)\cap L^2(\operatorname{div}\!=\!0;\Omega).
\end{align*}
It can be easily verified that these spaces are closed subspaces of $L^2(\Omega)^3$, $H(\operatorname{curl};\Omega)$, and $H_0(\operatorname{curl};\Omega)$, respectively. Therefore, $L^2(\operatorname{div}\!=\!0;\Omega)$ is a Hilbert space endowed with the $L^2(\Omega)^3$ norm, and both $H(\operatorname{curl},\operatorname{div}\!=\!0;\Omega)$ and $H_0(\operatorname{curl},\operatorname{div}\!=\!0;\Omega)$ are Hilbert spaces equipped with the standard norm of $H(\operatorname{curl};\Omega)$. Next we present two lemmas on the curl operator on the space of divergence-free functions.
\begin{lem}\label{lem:curlcoerc}
Let $\Omega\subset\R^3$ be a~simply connected and bounded Lipschitz domain with connected boundary. Then there exists some $c>0$, such that
\[\|\operatorname{curl}x\|_{L^2(\Omega)}\geq c\,\|x\|_{L^2(\Omega)}\quad \forall x\in H_0(\operatorname{curl},\operatorname{div}\!=\!0;\Omega).\]
\end{lem}
\begin{proof}
This result is a direct consequence of \cite[Cor.~3.51]{Monk2003}, where divergence-free functions with possibly nontrivial tangential traces are considered, and their $L^2$-norm is proven to be bounded from below by the sum of $L^2$-norm of the curl and the $L^2$-norm of the tangential trace.
\end{proof}
In the next step, we demonstrate that imposing the divergence-free constraint does not lead to a smaller space of admissible tangential traces.
\begin{lem}\label{lem:div0trace}
Let $\Omega \subset \R^3$ be a simply connected, bounded Lipschitz domain.
Then the restrictions of the trace operators $\gamma_\times$ and $\gamma_\pi$ to $H(\operatorname{curl},\operatorname{div}\!=\!0;\Omega)$ are surjective onto $V_\times(\partial\Omega)$ and $V_\pi(\partial\Omega)$, respectively.
\end{lem}
\begin{proof}
We restrict the proof to the case of $\gamma_\times$, as the argument for $\gamma_\pi$ is entirely analogous.

{\em Step~1:} We show that $\nabla\varphi\in H_0(\operatorname{curl};\Omega)$ for all $\varphi\in H^1_0(\Omega)$.\\
First, since $\operatorname{curl}\nabla \varphi=0$, we have that
$\nabla\varphi\in H(\operatorname{curl};\Omega)$ with
\begin{equation}
\|\nabla\varphi\|_{L^2(\Omega)}=\|\nabla\varphi\|_{H(\operatorname{curl};\Omega)}.\label{eq:curlnormgrad}
\end{equation}
By definition, there exists a~sequence $(\varphi_n)$ of compactly supported smooth functions converging in $H^1(\Omega)$ to $\varphi$.
Then, by using \eqref{eq:curlnormgrad}, we have that $(\nabla\varphi_n)$ converges in $H(\operatorname{curl};\Omega)$ to $\nabla\varphi$. Invoking that $\gamma_\times(\nabla\varphi_n)=0$ for all $n\in\N$, boundedness of the tangential trace operator now implies that $\gamma_\times(\nabla\varphi)=0$, and \eqref{eq:kergamma} yields that $\nabla\varphi\in H_0(\operatorname{curl};\Omega)$.\\
{\em Step~2:} We show the desired result. Let $f\in V_\times(\partial\Omega)$.
Then Theorem~\ref{thm:tangtraces}\,\ref{thm:tangtracesc} yields the existence of some
$x_1\in H(\operatorname{curl};\Omega)$ such that $\gamma_\times x_1=f$.
The Lax--Milgram lemma (see \cite[Thm.~6.2]{alt2016linear})  ensures the existence
of some $x_2\in H_0^1(\Omega)$ with
\begin{equation*}
\left\langle\nabla x_2,\nabla \psi\right\rangle_{L^2(\Omega)} = -\left\langle x_1,\nabla\psi\right\rangle_{L^2(\Omega)}
\quad\forall\, \psi\in H_0^1(\Omega).
\end{equation*}
We define $x \coloneqq  x_1 + \nabla x_2$.
Then the weak divergence of $x$ vanishes, since for all $\varphi\in H_0^1(\Omega)$, it holds that
\[\left\langle x,\nabla\varphi\right\rangle_{L^2(\Omega)}=\left\langle x_1,\nabla\varphi\right\rangle_{L^2(\Omega)}+\left\langle \nabla x_2,\nabla\varphi\right\rangle_{L^2(\Omega)}=0.\]
Moreover, since $\nabla x_2 \in H_0(\operatorname{curl};\Omega)$ by Step~1, we obtain that $x\in H(\operatorname{curl};\Omega)$ and
\[\gamma_\times(x)=\underbrace{\gamma_\times x_1}_{=f}+\underbrace{\gamma_\times x_2}_{=0}=f.\]
This proves that $\gamma_\times$ is surjective onto $V_\times(\partial\Omega)$.
\end{proof}

	\section{Skew-adjointness of $\calA$ from \eqref{eq:As_neu}}\label{app:supplementary}
	First, observe the decomposition
	\begin{align*}
		\calA = J F \calH , \qquad \text{where} \qquad F=\sbmat{}{I}{I}{}: \stack{\calX}{\calX_h} \to \calZ_h, \quad J=\sbmat{I}{}{}{-I}: \calZ_h \to \calZ_h.
	\end{align*}
	of $\calA$ defined in \eqref{eq:As_neu}.
	\begin{lem}\label{lem:As_skew}
		The operator $\calA$ from \eqref{eq:As_neu} is skew-adjoint, i.e $\calA=-\calA^*$.
	\end{lem}
	\begin{proof}
		We first show that $F\calH$ is a self-adjoint operator in $\calZ_h$. To this end, observe that
		\begin{align*}
			\langle F\calH \roundstack{z}{w}, \roundstack{x}{y} \rangle_{\calX_h \times \calX}&= \langle w,x \rangle_{\calX_h} + \langle S z, y \rangle_{\calX} =\langle  w , Sx \rangle_{\calX_h ,\calX_h^*} + \langle  Sz, y \rangle_{\calX_h^*, \calX_h}  \\
			&=\langle  w , Sx \rangle_{\calX} +  \langle  z, y \rangle_{\calX_h} = \langle  \roundstack{z}{w}, F\calH \roundstack{x}{y} \rangle_{\calX_h \times \calX}
		\end{align*}
		for all $\roundstack{z}{w}, \roundstack{x}{y} \in \dom F\calH= \dom \calA$ which shows $F\calH \subset (F\calH)^*$. Moreover, the operator
		\begin{align*}
			\begin{pmatrix}
				0 & S^{-1} \\ I & 0
			\end{pmatrix}: \calZ_h \to \calZ_h
		\end{align*}
		is a bounded inverse of $F\calH$ since
		\begin{align*}
			\begin{bmatrix}
				0 & S^{-1} \\ I & 0
			\end{bmatrix} \begin{bmatrix}
				0 & I \\ S & 0
			\end{bmatrix}  =\begin{bmatrix}
				0 & I \\ S & 0
			\end{bmatrix} \begin{bmatrix}
				0 & S^{-1} \\ I & 0
			\end{bmatrix} = I_{ \calZ_h}.
		\end{align*}
		This implies the self-adjointness by
		\begin{align*}
			(F\calH)^* = ((F\calH)^{-*})^{-1}=((F\calH)^{-1})^{-1}=F\calH.
		\end{align*}
		As a consequence
		\begin{align*}
			\calA^* = (J F\calH)^* =  (F\calH)^* J^* = JJF\calH J = J  \calA J = -\calA.
		\end{align*}
	\end{proof}

\bibliographystyle{abbrv}
\bibliography{main}

\newpage

\thispagestyle{empty}
\newgeometry{landscape, margin=2cm}

\begin{landscape}
    \section{Summary of involved spaces and operators}\label{sec:ourops}
	\vspace{0cm} \renewcommand{\arraystretch}{1.3}
    \centering
    \scalebox{0.8}{
\begin{tabular}[h]{|m{0.45cm}|m{4.5cm}|m{6.5cm}|m{3cm}|m{6.5cm}|}
\hline
& \textbf{Symbol} & \textbf{Meaning} & \textbf{Where defined?} & \textbf{Wave ex.} \\
\hline \hline
\multirow{6}{*}{\vspace{-2.2cm}\begin{turn}{90}\textbf{Spaces}\end{turn} }
\rule{0pt}{4mm} & $\calX$ & Hilbert space on which stiffness operator $S$ is defined & Sec.~\ref{sec:main}  & $L^2(\Omega)$ \\
\cline{2-5}
\rule{0pt}{6mm} & $\calX_h$ & $(\dom S^{1/2}, \|S^{1/2} \cdot \|_\calX)$ & Eq.~\eqref{eq:defX_h}  & $(H^1(\Omega), || \grad \cdot ||_{L^2(\Omega)})$ \\
\cline{2-5}\rule{0pt}{6mm}
 & $\calY$ & codomain of $A$ where $A^*A=S$ & Subsec.~\ref{subsec:abpot}  & $L^2(\Omega)^n$ \\
\cline{2-5}\rule{0pt}{6mm}
 & $\calY_1$ & $(\ran A, \| \cdot \|_\calY)$ & Eq.~\eqref{eq:Y_1}  & $(\grad H^1(\Omega),||\cdot ||_{L^2(\Omega)^n})$ \\
\cline{2-5}\rule{0pt}{6mm}
 & $\calX_M$ & $(\calX, \|M^{-1/2} \cdot \|_\calX)$ & Sec.~\ref{sec:main}  & $(L^2(\Omega), \|\rho^{-1/2} \cdot \|_{L^2(\Omega)})$ \\
\cline{2-5}\rule{0pt}{6mm}
 & $\calZ_h$ & product space $\stack{\calX_h}{\calX_M}$ & Eq.~\eqref{eq:Z_h}  & $\stack{H^1(\Omega)}{L^2(\Omega)}$ \\
\hline \hline
\multirow{9}{*}{\vspace{-5.15cm}\begin{turn}{90}\textbf{Operators}\end{turn} }
\rule{0pt}{6mm} & $S: \calX \supset \dom S \to \calX$ & stiffness operator & Sec.~\ref{sec:main}  & $-\Delta: L^2(\Omega) \supset \dom (\Delta) \to L^2(\Omega)$ \\
\cline{2-5}\rule{0pt}{6mm}
 & $M: \calX \to \calX$ & abstract mass density operator & Sec.~\ref{sec:main} & $x \mapsto \rho(\cdot) x: L^2(\Omega) \to L^2(\Omega)$ \\
\cline{2-5}\rule{0pt}{6mm}
 & $D: \calX \to \calX$ & dissipation operator & Sec.~\ref{sec:main} & $x \mapsto b(\cdot) x: L^2(\Omega) \to L^2(\Omega)$ \\
\cline{2-5}\rule{0pt}{6mm}
 & $A: \calX \supset \dom A \to \calY$ & factorization of $S=A^*A$ & Subsec.~\ref{subsec:abpot} &  $\grad: L^2(\Omega) \supset H^1(\Omega) \to L^2(\Omega)^n$  \\
\cline{2-5}\rule{0pt}{6mm}
 & $B: \calY \supset \dom B \to \calX$ & extension of $-A^*$ & Eq.~\eqref{eq:dual_pair_ass}  & $\divergence: L^2(\Omega)^n \supset H^{\divergence}(\Omega) \to L^2(\Omega)$ \\
\cline{2-5}\rule{0pt}{6mm}
 & $B_{\calY_1}: \calY_1 \supset \dom B_{\calY_1} \to \calX$ & restriction of $B$ to $\calY_1$ & Eq.~\eqref{eq:B_res}  & $\divergence:\grad H^1(\Omega) \supset H^{\divergence}(\Omega) \cap \grad H^1(\Omega) \to L^2(\Omega)$ \\
\cline{2-5}\rule{0pt}{8mm}
 & $A: \calX_h  \to \calY_1$ & unique bounded restriction of $A$ & Eq.~\eqref{eq:bounded_restriction}  & $\grad: H^1(\Omega) \to \grad H^1(\Omega)$ \\
\cline{2-5}\rule{0pt}{6mm}
 & $S: \calX_h  \to \calX_h^*$ & unique bounded restriction of $S$, $S=A'A$ & Subsec.~\ref{subsec:abpot}  & $-\Delta: H^1(\Omega) \to (H^1(\Omega))^*$ \\
\cline{2-5}\rule{0pt}{6mm}
 & $\calA: \calZ_h \supset \dom \calA \to \calZ_h$ & operator matrix corresponding to second order system & Eq.~\eqref{eq:As}  &  $\sbmat{0}{I}{\Delta}{-b}$  \\
\cline{2-5}\rule{0pt}{6mm}
 & $\calA_\mathrm{max}: \calZ_h \supset \dom \calA_\mathrm{max}\to \calZ_h$ & extension of $\calA$ in case $M=I$, $D=0$ & Eq.~\eqref{eq:defndomA}  & $\sbmat{0}{I}{\Delta}{0}$ \\
\cline{2-5}\rule{0pt}{6mm}
 & $\calB_1: \stack{\calY_1}{\calX} \supset \dom \calB_1 \to \stack{\calY_1}{\calX}$ & similarity transform of $\calA_\mathrm{max}$ & Eq.~\eqref{eq:defndomB_1}  & $\sbmat{0}{\grad}{\divergence}{0}$  \\
\hline
\end{tabular}}\hspace*{2cm}
\captionof{table}{Involved spaces and operators}\label{tab:ourops}
\end{landscape}
\restoregeometry
\newpage

\end{document}